\documentclass[journal,onecolumn]{IEEETran}

\usepackage[english]{babel}
\usepackage[utf8]{inputenc}

\usepackage{csquotes}
\usepackage{fancyhdr}

\usepackage{enumitem}
\usepackage{mathtools, amsmath, amssymb}
\usepackage{multirow}
\usepackage{picins}
\usepackage{cancel}
\usepackage{xcolor}
\usepackage{empheq}
\usepackage{hyperref}
\usepackage{diffcoeff}
\diffdef { p }
{
op-symbol = \partial ,
left-delim = \left .,
right-delim = \right | ,
subscr-nudge = +4 mu
}
\usepackage{caption}
\usepackage{subcaption}
\usepackage{pifont}



\newcommand{\R}{\mathbb{R}}
\newcommand{\B}{\mathbb{B}}
\newcommand{\N}{\mathbb{N}}

\newcommand{\abs}[1]{\left\lvert#1\right\rvert}

\newcommand*{\QEDB}{\hfill\ensuremath{\square}}%

\newcommand{\set}[1]{\left\{#1\right\}}

\newcommand*\tageq{\refstepcounter{equation}\tag{\theequation}}

\newcommand{\grad}[1]{\ensuremath{\text{grad}#1}}

\allowdisplaybreaks

\newtheorem{theorem}{Theorem}[section]
\newtheorem{corollary}[theorem]{Corollary}
\newtheorem{lemma}[theorem]{Lemma}

\newtheorem{example}[theorem]{Example}

\newtheorem{definition}{Definition}[section]{\bf}{\it}
\newtheorem{stdAssumption}{Standing Assumption}[section]
\newtheorem{remark}{Remark}[section]

\makeatletter
\DeclareRobustCommand{\qed}{%
  \ifmmode 
  \else \leavevmode\unskip\penalty9999 \hbox{}\nobreak\hfill
  \fi
  \quad\hbox{\qedsymbol}}
\newcommand{\openbox}{\leavevmode
  \hbox to.77778em{%
  \hfil\vrule
  \vbox to.675em{\hrule width.6em\vfil\hrule}%
  \vrule\hfil}}
\newcommand{\qedsymbol}{\openbox}

\newcommand{\proofname}{Proof}
\makeatother

\usepackage{graphicx}
\graphicspath{ {./} }

\begin{document}
\title{Model-Free Optimization on Smooth Compact Manifolds: Overcoming Topological Obstructions using Zeroth-Order Hybrid Dynamics} 

\author{Daniel E. Ochoa and Jorge I. Poveda  \vspace{-0.86cm}      
\thanks{The authors are with the Department of Electrical and Computer Engineering at the University of California San Diego, CA, USA. E-mail: {\tt dochoatamayo@ucsd.edu}.}}
\maketitle


\begin{abstract}
Smooth autonomous dynamical systems modeled by ordinary differential equations (ODEs) cannot robustly and globally stabilize a point in  compact, boundaryless manifolds. This obstruction, which is topological in nature, implies that traditional smooth optimization dynamics are not able to robustly solve global optimization problems in such spaces. In turn, \emph{model-free} optimization algorithms, which usually inherit their stability and convergence properties from their model-based counterparts, might also suffer from similar topological obstructions. For example, this is the case in zeroth-order methods and perturbation-based techniques, where gradients and Hessians are usually estimated in real time via measurements or evaluations of the cost function. To address this issue, in this paper we introduce a class of \emph{hybrid model-free optimization dynamics} that combine continuous-time and discrete-time feedback to overcome the obstructions that emerge in traditional ODE-based optimization algorithms evolving on smooth compact manifolds. In particular, we introduce a hybrid controller that switches between different model-free feedback-laws obtained by applying suitable exploratory \emph{geodesic dithers} to a family of synergistic diffeomorphisms adapted to the cost function that defines the optimization problem. The geodesic dithers enable a suitable exploration of the manifold while preserving its forward invariance, a property that is important for many practical applications where physics-based constraints limit the feasible trajectories of the system. The hybrid controller exploits the information obtained from the geodesic dithers to achieve robust global practical stability of the set of minimizers of the cost. The proposed method is of model-free nature since it only requires \emph{measurements} or evaluations of the cost function. Numerical results are presented to illustrate the main ideas and advantages of the method.
\end{abstract}

\section{Introduction}
\label{sec:intro}
This paper studies algorithms for the \emph{global} solution of optimization problems of the form
\begin{align*}
    \inf~ \phi(z)\quad\text{ subject to }\quad z\in M,\tageq{\label{statement:optimization}}
\end{align*}
where $(M,g)$ is an $n$-dimensional Riemannian manifold to be formally defined in Section \ref{sec:preliminaries}, and $\phi$ is a smooth cost function. The mathematical forms of $\phi$ and its derivatives are assumed to be \emph{unknown}, and we only assume that $\phi$ is available via \emph{measurements} or \emph{evaluations}. Optimization problems on manifolds, of the form \eqref{statement:optimization}, emerge in several interesting applications, ranging from statistics \cite{Freifeld:Thesis:2013} to robotics \cite{watterson2020trajectory}, aerospace engineering \cite{hauser2002projection}, power systems \cite{absil2009optimization}, and quantum control \cite{Bassam}. Among the most simple and successful algorithms for optimization, gradient-descent methods have been studied in the context of manifolds since at least the end of last century \cite{gabay1982minimizing}. Such methods have also recently garnered considerable interest due to their potential in estimation, machine learning, and data science pipelines \cite{bottou2010large}.  In the context of dynamical systems described by ordinary differential equations (ODEs), real-time optimization problems defined on manifolds are typical in robotics, mechanical systems, and aerospace control problems evolving under kinematic constraints, e.g., unicycle \cite[Sec. 2.2]{sontag1999stability}, or obstacle-occluded spaces  \cite{poveda2021robust}. In such problems, the restriction to evolve on a smooth manifold limits the feasible directions that any onboard optimization algorithm can exploit in real time. For comprehensive introductions to ODE-based optimization algorithms on manifolds, we refer the reader to \cite{helmke2012optimization,absil2009optimization}.

One of the key challenges in the solution of optimization problems on smooth (boundaryless) compact manifolds is given by the fact that in such spaces a point cannot be \emph{robustly globally} asymptotically stabilized by using continuous feedback in ODEs \cite[Thm. 1]{bhat2000topological}, a result that actually extends to Lie groups and non-contractible spaces in general, see \cite[Thm. 21]{sontag2013mathematical}. This well-known result implies that standard gradient-descent or Newton-like flows cannot achieve robust \emph{global} convergence to the minimizer of a continuously differentiable cost function for every type of smooth compact manifold. The reason behind this incompatibility lies in the fundamental difference between the topological nature of the basin of attraction of a point under continuous dynamics, and that of a compact boundaryless manifold \cite[Thm. 21]{sontag2013mathematical}: the former is contractible while the latter is not.

Existing results in the literature circumvent this issue by using notions of asymptotic stability that disregard measure-zero sets containing critical points of the cost function that are not solutions of the optimization problem under study \cite{angeli2001almost,angeli2004almost,efimov2012global}, e.g., local maximizers, saddle points, etc. However, algorithms with \emph{almost} global convergence certificates have been shown to be susceptible to arbitrarily small (adversarial) disturbances, under which the set of initial conditions from which convergence is not achieved is not of measure zero anymore, see \cite[Thm. 6.5]{sanfelice2007thesis}, \cite{mayhew2011synergistic}, \cite{sontag1999stability},\cite[Ex. 1]{poveda2021robust}, and \cite[Ex. 5.2]{sanfelice2006robust} for specific examples.

Alternatively, there have been developments that confront the obstruction by using time-varying \cite{coron1992global,sepulchre1992some}, or discontinuous feedback \cite{ancona1999patchy,malisoff2006global}, finding success in recovering global stability certificates under some mild assumptions. However, as shown in \cite[Cor. 21]{mayhew2011topological}, these results only circumvent the issue when the optimization dynamics operate in nominal conditions. Namely, when disturbances are unavoidable (and adversarial), \emph{robust} and global stabilization of a point in compact boundaryless manifolds cannot be achieved by merely using discontinuous or time-varying feedback strategies. Since guaranteeing positive margins of robustness with respect to additive disturbances is of fundamental importance for practical implementations where noisy measurements are unavoidable, in \cite{mayhew2010hybrid} the authors introduced a hybrid controller able to overcome the obstructions. The proposed strategy makes use of hybrid dynamical systems theory \cite{bookHDS} to prescribe feedback control dynamics that synergistically switch between different continuous vector fields generated from a family of potential functions. Recent works have employed the synergistic framework to solve attitude stabilization problems in aerospace application \cite{mayhew2011quaternion}, stabilization by hybrid backstepping \cite{mayhew2011synergistic}, stabilization in $\text{SO}(3)$ \cite{berkane2017hybrid}, and for the robust stabilization of trajectories in multi-rotor aerial vehicles \cite{casau2019robust}. However, since these works address stabilization problems, where the point to be stabilized is known \emph{a priori}, they cannot be directly used for the solution of optimization problems where the set of optimizers is unknown. 

On the other hand, in optimization problems where the cost function is unknown and only available via measurements or evaluations, model-free optimization algorithms of \emph{zeroth-order} are required. While the literature of continuous-time zeroth-order dynamics, also known as extremum seeking algorithms, is quite rich \cite{krstic2000stability,tan2006non,durr2014extremum,poveda2017framework}, most of the results in the literature for model-free optimization on smooth compact manifolds have been developed for smooth dynamical systems that seek to emulate the behavior of smooth gradient flows \cite{durr2014extremum,taringoo2018optimization,poveda2015shahshahani}. Since in those settings the ``practical'' asymptotic stability properties of the model-free dynamics are completely inherited from the asymptotic stability properties of the smooth model-based algorithms, the obstructions to robust global optimization extend to the model-free counterparts. Additionally, existing results that can achieve global convergence via switching control \cite{strizic2017hybrid} are not able to guarantee the forward invariance of the manifold by introducing dither signals that do not evolve on its tangent space, a requirement that is relevant for applications where the evolution on manifolds is enforced by physical constraints, or in problems where the cost function is defined only on the manifold. 

With this in mind, the main contribution of this paper is the introduction of a novel class of model-free optimization algorithms for problems defined on compact boundaryless connected Riemannian manifolds that impose topological obstructions to global optimization using smooth dynamics. In particular, we prescribe a family of hybrid model-free controllers that switch between different zeroth-order feedback laws obtained by applying suitable exploratory geodesic dithers to a class of synergistic family of diffeomorphisms adapted to the cost function $\phi$. Under such dynamics, we present theoretical \emph{robust global} asymptotic stability certificates for a neighborhood of the set of minimizers (i.e., practical stability), by employing state-of-the-art results on hybrid model-free dynamics \cite{bookHDS,poveda2017framework,Poveda2021RobustHZ}. The proposed approach performs \emph{exploration} of the state space via geodesic dithers to obtain -on average- suitable evolution directions, and simultaneous \emph{exploitation} of these directions via hybrid controllers that guarantee \emph{robust global practical} asymptotic stability of the set of minimizers of the cost. Since our algorithms are hybrid and do not necessarily have the uniqueness of solutions property, we leverage averaging techniques for non-smooth and hybrid systems that guarantee the closeness of solutions property (on compact time domains) between the original hybrid dynamics and some solution of the average hybrid  dynamics. This closeness of solutions property, studied for nonsmooth and hybrid systems in \cite[Thm. 1]{wang2012analysis} and \cite[Prop. 6]{Poveda2021RobustHZ} relies on sequential compactness results that use the graphical distance metric \cite{goebel2006solutions} instead of the standard Arzela-Ascoli theorem used for ODEs, which in general is not applicable to hybrid systems. Compared to previous approaches for model-free optimization on manifolds, e.g. \cite{taringoo2018optimization}, \cite{durr2014extremum} and \cite{suttner2022extremum}, our convergence results are global rather than local or almost global. 

Compared to existing switching algorithms \cite{strizic2017hybrid}, our model-free dynamics are designed to evolve on the manifold and preserve its invariance via geodesic dithering. The results presented in this paper are also applicable to a larger class of manifolds and optimization problems. Our results also provide an alternative approach to the solution of model-free optimization and extremum seeking problems with multiple critical points, typical in non-convex settings, a problem that has also been recently studied in \cite{krstic2022} and \cite{marconi2022} using other techniques.

The rest of this paper is organized as follows. Section \ref{sec:preliminaries} presents preliminary notions on hybrid systems and Riemannian geometry. Section \ref{sec:main} presents the main results, including the general hybrid model-free optimization dynamics and two specific examples of algorithms synthesized for two different applications. We continue by presenting the analysis and proofs in Section \ref{sec_analysis}. Finally, we end by presenting the conclusions in Section \ref{sec:conclusions}.

\section{Preliminaries}
\label{sec:preliminaries}

In this section, we introduce the notation used in the paper, as well as some preliminary notions on Riemannian geometry and hybrid dynamical system's theory.
\subsection{Notation}

Given a compact set $\mathcal{A}\subset N$ in a metric space $N$, with metric $d:N\times N\to \R_{\ge 0}$, and an element $z\in N$, we use $|z|_{\mathcal{A}}\coloneqq \min_{s\in\mathcal{A}}d(z,s)$ to denote the minimum distance of $z$ to $\mathcal{A}$.  We use $\mathbb{S}^n\coloneqq \{z\in\R^{n+1}:\sum_{i=1}^{n+1} z_i^2=1\}$ to denote the $n$-th dimensional sphere, with $\mathbb{S}^1$ representing the unit circle in $\R^2$. Following suit, we use $\mathbb{T}^n=\mathbb{S}^1\times\ldots\times\mathbb{S}^1$ to denote the $n^{th}$ Cartesian product of $\mathbb{S}^1$. We also use $r\mathbb{B}$ to denote a closed ball in the Euclidean space, of radius $r>0$, and centered at the origin. We use $I_n\in\R^{n\times n}$ to denote the identity matrix. A function $\beta:\R_{\geq0}\times\R_{\geq0}\to\R_{\geq0}$ is said to be of class $\mathcal{K}\mathcal{L}$ if it is non-decreasing in its first argument, non-increasing in its second argument, $\lim_{r\to0^+}\beta(r,s)=0$ for each $s\in\R_{\geq0}$, and  $\lim_{s\to\infty}\beta(r,s)=0$ for each $r\in\R_{\geq0}$. We use $\pi_A : A\times B\to A$ to denote the natural projection from the product space $A\times B$ to $A$, and $\text{gph}~J$ to denote the graph of a mapping. The Kronecker delta is denoted with $\delta^i_j$. Throughout the paper we make use of Einstein summation convention \cite[Ch. 1, pp. 18]{lee2013smooth}.
\subsection{Riemannian Manifolds}

We introduce the main differential geometric concepts used throughout the paper. For a more thorough discussion of these topics, we refer the reader to \cite{lee2013smooth} and \cite{lee2018riemann}. The concept of smooth manifold will play an important role in this paper:

\textbf{Smooth manifolds:} An $n$-dimensional manifold is a second-countable Hausdorff topological space that is locally Euclidean of dimension $n$. A coordinate chart for $M$ is a pair $(U, \varphi)$ where $U\subset M$ is an open set and $\varphi:U\to \hat{U}\subset\R^n$ is a homeomorphism. Two coordinate charts $(U,\varphi)$ and $(V, \psi)$ are said to be smoothly compatible if the transitions maps $\psi\circ \varphi^{-1}$ and $\varphi\circ \psi^{-1}$ are diffeomorphisms. A smooth structure on $M$ is a maximal collection of coordinate charts for which any two charts are smoothly compatible; a smooth coordinate chart is any chart that belongs to a smooth structure. Then, a \emph{smooth manifold} is a manifold endowed with a particular smooth structure. Given a smooth manifold $M$, the set of all smooth real-valued functions $f:M\to \R$ is denoted by $C^\infty(M)$.

Dynamical systems that evolve on smooth manifolds are characterized by vector fields that lie on their tangent space. Below, we formalize these concepts for a general manifold $M$:

\textbf{Tangent space and Vector Fields}: For every point $z\in M$, a tangent vector at $z$ is a linear map $v:C^\infty(M)\to \R$ that satisfies the product rule $v(fh) = f(z)\cdot v(h) + h(z)\cdot v(f)$, for $f,h\in C^\infty(M)$. The set of all tangent vectors at $z$ is denoted by $T_zM$, and is called the \emph{tangent space} of $M$ at $z$. For a given smooth cordinate chart $(U,\varphi)$ and by writing $\varphi = (z^1, \ldots, z^n)$, with $z^i:U\to\R$ the $i$-th coordinate map, the \emph{coordinate vectors} $\set{\diffp{}{z^1}[z],\ldots, \diffp{}{z^n}[z]}$ are defined as:
\begin{align*}
    \diffp{}{z^i}[z]f =  \diffp{}{z^i}[\varphi(z)](f\circ \varphi^{-1}),\tageq{\label{def:coordinate:vectors}}
\end{align*}
for an arbitrary $f\in C^\infty(M)$. The coordinate vectors form a basis for the $n$-dimensional vector space $T_zM$.\\
The tangent bundle of $M$ is defined to be the disjoint union of the tangent spaces at all points in the manifold,  $TM \coloneqq \bigsqcup_{z\in M}T_zM$. It has a natural topology and smooth structure inherited from $M$ that make it into a $2n$-dimensional smooth manifold. A \emph{smooth vector field} is a smooth map $X:M\to TM$ satisfying $X(z)\in T_z M$ for all $z\in M$, where $TM$ has the aforementioned smooth structure inherited from $M$. We use $\mathfrak{X}(M)$ to denote the set of all smooth vector fields on $M$.\\

A local frame for $M$ is an ordered tuple of vector fields $\left(X_1, \ldots, X_n\right)$  defined on an open set $U\subset M$ that is linearly independent and spans $T_zM$ at each $z\in M$. For any smooth chart $(U, (x^i))$, the assignment $z\to \diffp{}{z^i}[z]$ determines a smooth vector field on $U$, called the $i$-th coordinate vector field. The set of $n$ smooth coordinate vector fields constitutes a local frame on $U$, called the \emph{coordinate frame}.\\

For every $z\in M$, the cotangent space at $z$, denoted by $T_z^*M$, is defined to be the dual vector space to $T_zM$. Elements of $T_z^*M$ are called tangent covectors at $z$. Analogous objects to the tangent bundle, vector fields, and local frames can be defined with tangent covectors instead; they receive the names of cotangent bundle $T^*M$, covector fields, and local coframes, respectively. The differential of a function $f\in C^\infty(M)$, denoted by $df:TM\to \R$, is a covector field  defined pointwise by:
\begin{equation}\label{sec:prelim:differential}
    df_z(v)=vf,\qquad\forall~ v\in T_zM.
\end{equation}
For any smooth chart $(U,(z^i))$ the differentials, denoted $(dz^1,\ldots, dz^n)$ make up a local coframe on $U$ satisfying $dz^i\left(\diffp{}{z^j}\right)= \delta^i_j$ .

In this paper, we will restrict our attention to manifolds with a Riemannian structure:

\textbf{Riemannian Manifolds:} An $n$-dimensional \textit{Riemannian manifold} is a pair $(M, g)$, where $M$ is an $n$-dimensional smooth manifold, and $g$ is a Riemannian metric whose value at each point $z\in M$ is an inner product defined on the tangent space $T_zM$. Given smooth local coordinates $(z^1,\ldots, z^n)$ on an open subset $U\subseteq M$, the Riemannian metric can be written locally as:
$g = g_{ij}dz^idz^j,$
where $dz^idz^j$ represents the symmetric tensor product between the coordinate coframe elements $dz^i$ and $dz^j$, $i,j\in\set{1,\ldots,n}$, and $(g_{ij})$ is a symmetric positive definite matrix of smooth functions $g_{ij}\in C^\infty(M)$, characterized by $g_{ij}(z) = g\left(\diffp{}{z^i}[z],\diffp{}{z^j}[z]\right)$.

\textbf{Critical points of scalar functions:} The sets of \textit{critical points and values} of $f\in C^\infty(M)$ are respectively defined by
\begin{align}
    \text{Crit}f &\coloneqq \set{z\in M~:~df_z=0},\label{sec:prelim:criticalPoints}\\
    \text{Val}f &\coloneqq \set{a \in \R~:~a=f(z),~z\in \text{Crit}f}\label{sec:prelim:criticalValues}.
\end{align}
The gradient of $f$ is a continuous map denoted by $\grad{f}:M\to TM$, and characterized by 
\begin{equation}
    df_z(v) = g\left(\grad{f}|_z,v\right),~~ \text{for all }z\in M, v\in T_zM,\tageq{\label{prelim:gradient}}
\end{equation}
where $\grad{f}|_z\in T_zM$ represents the value of the gradient of $f$ at $z$. With the definition of the gradient in \eqref{prelim:gradient}, we can assess whether a point $z$ is a critical point of $f$ or not, by checking if the gradient of $f$ vanishes identically at that point.

To guarantee suitable exploration of $M$, while preserving its forward invariance, we will work with algorithms that implement geodesic dithers:

\textbf{Geodesics:} Geodesics are defined as curves $\gamma:[a,b]\to M$ on a Riemannian manifold, satisfying
\begin{align*}
    \nabla_{\dot{\gamma}(t)}\dot{\gamma}(t) = 0,\tageq{\label{def:geodesic}}
\end{align*}
where $\nabla:\mathfrak{X}(M)\times \mathfrak{X}(M)\to \mathfrak{X}(M)$ is the Levi-Civita connection \cite[Ch. 5]{lee2018riemann}. By expressing \eqref{def:geodesic} in coordinates, the component functions $\gamma_i$ of a geodesic $\gamma$ can be written as the solution of the  differential equations:
\begin{align*}
    &\ddot{\gamma}_i(s) + \Gamma_{j,k}^i\dot{\gamma}_j(s)\dot{\gamma}_k(s) =0,\quad i\in\set{1,\ldots,n},\\
    &\Gamma_{ij}^k = \frac{1}{2} g^{il}\left(\diffp{g_{jl}}{z^k}+\diffp{g_{kl}}{z^j}-\diffp{g_{jk}}{z^l}\right),
\end{align*}
where $g^{ij}\coloneqq (g_{ij})^{-1}$, and $g_{ij}$ are the local components of the Riemannian metric $g$. To generate the dither signals used by the model-free optimization algorithms considered in this paper, we use the \emph{restricted exponential map} $\text{exp}_z: T_zM\to M$, defined by $\text{exp}_z(v)= \gamma_v(1)$, where $\gamma_v$ is the unique maximal geodesic satisfying $\gamma(0)=z$ and $\dot{\gamma}(0)=v$.

Throughout the paper we make use of the following standing assumption.

\begin{stdAssumption}\label{assumption:topologicalPropertiesM}
The Riemannian manifold $(M,g)$ is compact, boundaryless, and connected. 
\end{stdAssumption}
\noindent In particular, Assumption \ref{assumption:topologicalPropertiesM} guarantees the existence of a path between any two points in the manifold \cite[Prop 2.50]{lee2018riemann}, which facilitates the definition of a notion of distance.

\textbf{Riemannian Distance:} The Riemannian distance, denoted by $d_g(z_1,z_2)$ is defined to be the infimum of the lengths of all admissible curves between a pair of points in the manifold \cite[Ch 2.]{lee2013smooth}. Formally, the Riemannian distance $d_g:M\times M\to\R_{\ge0}$ is defined by:
\begin{align*}
    d(z_1,z_2) \coloneqq \inf_{\gamma \in \mathcal{P}(z_1,z_2)} \int_{t_1}^{t_2} \sqrt{g\left(\dot{\gamma}(t),\dot{\gamma}(t)\right)}dt,\tageq{\label{riemannian:distance}}
\end{align*}
where $\mathcal{P}(z_1,z_2)$ represents the set of all admissible curves connecting $z_1$ and $z_2$, and $t_1,t_2\in\R$ are such that $\gamma(t_1)=z_1$ and $\gamma(t_2)=z_2$ for $\gamma\in \mathcal{P}(z_1,z_2)$.
\subsection{Hybrid Dynamical Systems and Stability Notions}

Since smooth autonomous ODEs cannot render a point robustly globally asymptotically stable on smooth compact manifolds \cite[Thm. 1]{bhat2000topological}, in this paper we consider model-free optimization algorithms modeled as hybrid dynamical systems (HDS) \cite{bookHDS} of the form:

\begin{subequations}\label{HDS}
\begin{align}
&x\in C,~~~~~~~\dot{x}= F(x)\label{flows00}\\
&x\in D,~~~~~x^+\in G(x),\label{jumps00}
\end{align}
\end{subequations}

where $x\in M$ is the state, $F:M\to TM$ is called the flow map, and $G:M\rightrightarrows M$ is a set-valued map called the jump map. The sets $C$ and $D$, called the flow set and the jump set, respectively, characterize the points in $M$ where the system can flow or jump via equations \eqref{flows00} or \eqref{jumps00}, respectively. Then, the HDS $\mathcal{H}$ is defined as the tuple $\mathcal{H}:=\{C, F, D, G\}$. Systems of the form \eqref{HDS} generalize purely continuous-time systems and purely discrete-time systems. Namely, continuous-time dynamical systems (e.g., ODEs) can be seen as a HDS of the form \eqref{HDS} with $D=\emptyset$, while discrete-time dynamical systems (e.g. recursions) correspond to the case when $C=\emptyset$.
Solutions to HDS of the form \eqref{HDS} are defined on hybrid time domains, i.e., they are parameterized by both a continuous-time index $t\in\mathbb{R}_{\geq0}$, and a discrete-time index $j\in\mathbb{Z}_{\geq0}$. Consequently, the notation $\dot{x}$ in \eqref{flows00} represents the derivative of $x$ with respect to time $t$, i.e., $\frac{dx(t,j)}{dt}$; and $x^+$ in \eqref{jumps00} represents the value of $x$ after an instantaneous jump, i.e., $x(t,j+1)$. For a precise definition of hybrid time domains and solutions to HDS of the form \eqref{HDS} we refer the reader to \cite[Ch.2]{bookHDS}. A HDS $\mathcal{H}$ is said to be well-posed if $C$ and $D$ are closed sets, $C\subset\text{dom}(F)$ and $D\subset\text{dom}(G)$, $F$ is continuous in $C$, and $G$ is outer-semicontinuous \cite[Def. 5.9]{bookHDS} and locally bounded \cite[Def. 5.14]{bookHDS} relative to $D$.\\

\textbf{Stability notions:} By endowing the manifold with the distance function $d_g$, $M$ constitutes a metric space \cite[Thm 2.55]{lee2018riemann}. Accordingly, we can use stability notions analogous to those studied in the traditional Euclidean space.

\begin{definition}\normalfont
The compact set $\mathcal{A}\subset C\cup D$ is said to be \emph{uniformly globally asymptotically stable} (UGAS) for system \eqref{HDS} if $\exists$ $\beta\in\mathcal{K}\mathcal{L}$ such that every solution $x$ satisfies:  
\begin{equation}\label{KLbound}
|x(t,j)|_{\mathcal{A}}\leq \beta(|x(0,0)|_{\mathcal{A}},t+j),
\end{equation}
$\forall~(t,j)\in\text{dom}(x)$, where $|z|_{\mathcal{A}}= \min_{s\in\mathcal{A}}d_g(z,s)$.
\QEDB
\end{definition}
\noindent We will also consider $\varepsilon$-parameterized HDS $\mathcal{H}_{\varepsilon}$:
\begin{align*}
x\in C_{\varepsilon},~~\dot{x}=F_{\varepsilon}(x),~~~\text{and}~~~x\in D_{\varepsilon},~x^+\in G_{\varepsilon}(x),
\end{align*}
where $\varepsilon>0$. For these systems we will study \emph{global practical stability} properties as $\varepsilon\to0^+$.

\begin{definition}\label{definitionSGPAS}\normalfont
The compact set $\mathcal{A}\subset C\cup D$ is said to be \emph{Globally Practically Asymptotically Stable} (GP-AS) as $\varepsilon\to0^+$ for system \eqref{HDS} if $\exists$ $\beta\in\mathcal{K}\mathcal{L}$ such that for each $\nu{>}0$ there exists $\varepsilon^*>0$ such that for all $\varepsilon\in(0,\varepsilon^*)$ and $x(0,0)\in M$, every solution of $\mathcal{H}_{\varepsilon}$ satisfies  
\begin{equation}\label{SGPASbound}
|x(t,j)|_{\mathcal{A}}\leq \beta(|x(0,0)|_{\mathcal{A}},t+j)+\nu,
\end{equation}
$\forall~(t,j)\in\text{dom}(x)$. ~\hfill\QEDB
\end{definition}
The notion of GP-AS can be extended to systems that depend on two parameters $\varepsilon=(\varepsilon_1,\varepsilon_2)$. In this case, we say that $\mathcal{A}$ is GP-AS as $(\varepsilon_{2},\varepsilon_{1})\to0^+$ where the parameters are tuned in order starting from $\varepsilon_1$. 
%

\section{Main Results}
\label{sec:main}
Approaches for optimization in Euclidean spaces with global convergence certificates usually rely on convexity properties of the cost function $\phi$. For Riemannian manifolds, convexity is characterized along geodesics. However, this concept turns out to be of little utility under the light of Assumption \ref{assumption:topologicalPropertiesM}. Indeed, compact Riemannian manifolds do not admit nonconstant geodesically convex  functions \cite[Cor. 2.5]{udriste2013convex}. In view of the limitations imposed by assuming convexity of functions in compact Riemannian manifolds, in this paper we instead make use of the following regularity assumption on $\phi$.

\begin{stdAssumption}\label{assumption:phi:critical}
The cost function $\phi$ has a finite amount of critical values, i.e., there exists $l\in\N$ such that that $\text{Val}~\phi = \set{\phi_1,\phi_2,\ldots,\phi_l}$, where $$\underline{\phi}\coloneqq \phi_1<\phi_2\leq \phi_i\leq \phi_l\eqqcolon\overline{\phi},$$ and $\underline{\phi}\leq \phi(z)\leq \overline{\phi}$ for all $z\in M$. Moreover, $\phi$ has a unique minimizer.
\end{stdAssumption}
Let $\mathcal{A} \coloneqq \set{z\in \text{Crit}~\phi~:~ \phi(z)=\underline{\phi}}$ and $\mathcal{B}\coloneqq\text{Crit}~\phi\setminus\mathcal{A}$, represent the minimizer of $\phi$ and the set of critical points of $\phi$ that are not minimizers, respectively. Since $M$ is compact, the set $\mathcal{A}$ is also compact. Note that Assumption \ref{assumption:phi:critical} does not rule out functions with multiple first-order critical points where the gradient of $\phi$ vanish. Indeed, in our problem setup, $\mathcal{B}$ is not empty since, by Morse theory, there exists at least two critical points for scalar-valued functions on compact boundaryless manifolds \cite{morse1934calculus}. Such critical points will correspond to equilibria of traditional gradient flows, making them prone to arbitrarily small disturbances able to keep the trajectories of the flows away from the minimizers of $\phi$. This robustness issue, thoroughly discussed in \cite[Thm. 6.5]{sanfelice2007thesis}, \cite{mayhew2011synergistic}, \cite{sontag1999stability},\cite[Ex. 1]{poveda2021robust}, and \cite[Ex. 5.2]{sanfelice2006robust}, and further illustrated later in Section \ref{secrobustness} via numerical examples, is one of the main motivations for the development of switching algorithms that effectively rule out undesirable critical points in the system. In our case, we will design the hybrid algorithms to rule out undesirable critical points in the \emph{average dynamics} of the controllers, and we will show that this property will be retained by the original hybrid dynamics.
\begin{remark}\normalfont
For the particular case in which $\phi$ is a Morse function \cite{morse1934calculus}, Assumption \ref{assumption:phi:critical} is automatically satisfied. Moreover, since the set of Morse functions is known to be an open dense set in the space of differentiable functions, and provided some approximation error is admissible in our solutions, we can alternatively dispense of Standing Assumption \ref{assumption:phi:critical} by treating a surrogate approximate optimization problem to \eqref{statement:optimization}, whose solution is the minimizer of a Morse function that is sufficiently close to $\phi$. \QEDB 
\end{remark}
\subsection{Description of the Proposed Algorithms}
\begin{figure*}[t!]
\centering
\subfloat[Block diagram of the Zeroth-Order Hybrid Dynamics $\mathcal{H}_0$.]{%
  \includegraphics[width=0.49\textwidth]{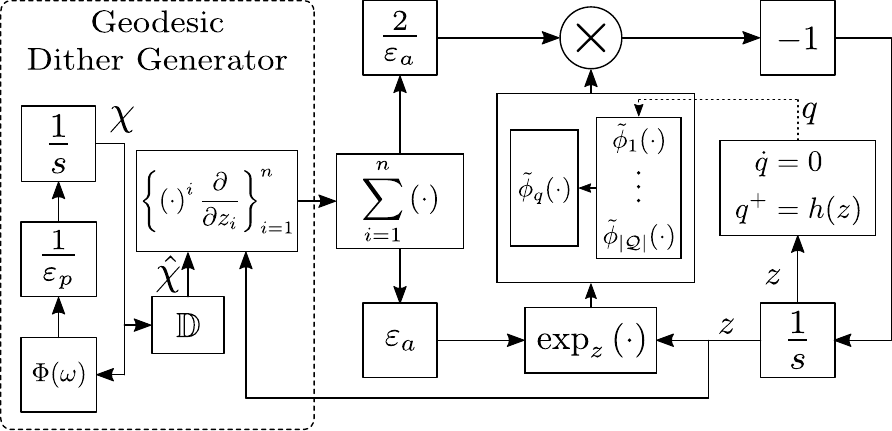}%
  \label{fig:schemeGeodesicES}%
}\qquad
\subfloat[Example of a trajectory $z$ and its average $z^A$.]{%
  \includegraphics[width=0.43\textwidth]{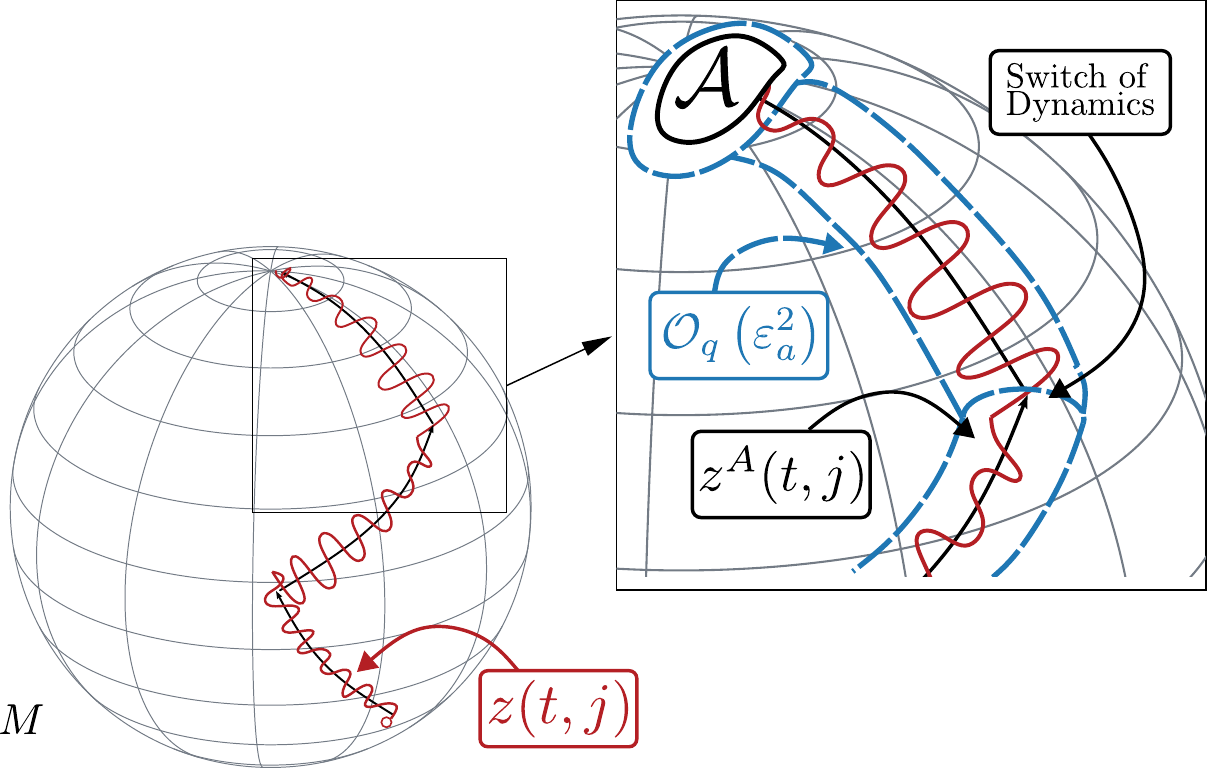}%
  \label{fig:cartoon}%
}
\caption{Left: Block diagram of the proposed hybrid zeroth-order dynamics with geodesic dithering. Right: Cartoon of the trajectories of the system evolving on a manifold M. The behavior of the oscillating trajectories is approximated by the trajectories of the average, first-order hybrid system $\mathcal{H}_1$.}
\label{figure1}
\end{figure*}

The left plot of Figure \ref{figure1} shows a block diagram describing the proposed controllers for the solution of problem \eqref{statement:optimization}. Before analyzing the mathematical properties of this system, we briefly describe the main ideas behind the controllers:
\begin{itemize}
\item  We make use of a set of dynamic oscillators to generate oscillating amplitudes whose dynamics are captured by the state $\chi$. The amplitudes are then used to drive dithering paths that extend along geodesics of the manifold $M$. These dithers will be used for the purpose of \emph{local exploration}.
\item The geodesic dithers, in conjunction with measurements of the cost $\phi$, are employed to generate a family of vector fields $\{\hat{f}_q(z,\chi)\}_{q\in \mathcal{Q}}$, given by
\begin{equation*}
\hat{f}_{q}\left(z,\chi\right) \coloneqq    \frac{2}{\varepsilon_a}\left(\tilde{\phi}_q\right)\left(\text{exp}_z\left(\varepsilon_a\hat{\chi}^i\diffp{}{z_i}\right)\right)\hat{\chi}^j\diffp{}{z_j}.
\end{equation*}
where $\varepsilon_a>0$ is a tunable gain. These vector fields, explained in detailed below, will be used for the purpose of \emph{exploitation} using zeroth-order dynamics.

\item In order to generate the family of vector fields $\{\hat{f}_q(z,\chi)\}_{q\in \mathcal{Q}}$, we use a family of $|\mathcal{Q}|$ diffeomorphisms that warp the manifold $M$, shifting the points that are not in a neighborhodod of the minimizers of $\phi$, and generating a family of warped cost functions $\{\tilde{\phi}_q\}_{q\in \mathcal{Q}}$. 

\item By \emph{partitioning} the manifold $M$ in a suitable way, we can implement each of the vector fields $\hat{f}_{q}\left(z,\chi\right)$ in ``safe zones'', described by regions where their \emph{average dynamics} (computed along the trajectories of the dither) are sufficiently far away from the critical points that are not optimal for the warped cost $\tilde{\phi}_q$.

\item As we increase the frequency of the dithers, the trajectories induced by the switching vector fields $\hat{f}_{q}\left(z,\chi\right)$ will approximate on compact time domains the trajectories of a class of switched gradient flows (i.e., first-order dynamics) able to achieve \emph{robust global} asymptotic stability of $\mathcal{A}$ on $M$.

\item By using averaging theory for well-posed hybrid inclusions, we will establish that the \emph{hybrid zeroth-order dynamics} will retain the global asymptotic stability properties of the average hybrid first-order methods, in a practical sense, i.e., convergence will be global but towards a neighborhood of $\mathcal{A}$.
\end{itemize}

The above ideas suggest that the proposed algorithms are similar in spirit to other hybrid controllers studied in the context of robust global \emph{stabilization} problems \cite{mayhew2010hybrid}, \cite[Ch. 7]{sanfelice2020hybrid}. However, the algorithms studied in this paper do not exactly fit the setting of synergistic hybrid control, since the family $\{\hat{f}_q(z,\chi)\}_{q\in \mathcal{Q}}$ does not correspond to the gradients of synergistic Lyapunov functions. In fact, unlike standard stabilization problems, the main challenges in problem \eqref{statement:optimization} are that set $\mathcal{A}$ and the function $\phi$ are unknown. Therefore, to implement the model-free hybrid dynamics we need to characterize the family of cost functions $\phi$ and smooth manifolds $(M,g)$ that admit suitable partitions and deformations to generate feasible switching rules that induce global stability to $\mathcal{A}$, in a model-free way. This characterization is the main subject of the next section.
\subsection{Stability, Convergence, and Robustness Results}
The closed-loop system describing the zeroth-order hybrid dynamics shown in Figure \ref{fig:schemeGeodesicES} makes use of three main states $y \coloneqq (z,~q,~\chi)\in M\times \mathcal{Q}\times \mathbb{T}^{n}$, where $z$ is an internal auxiliary state, $q\in \mathcal{Q}:=\{1,2,\ldots,N\}$, $N\in\mathbb{Z}_{\geq2}$, is a logic decision variable, and $\chi$ is the state of the oscillator. The data of this hybrid system is denoted as:
\begin{equation}\label{hybridzeroorder}
\mathcal{H}_{0} = \{C_{0},~F_{0},~ D_{0},~G_{0}\}.
\end{equation}
The continuous-time dynamics of $\mathcal{H}_{0}$ are given by the constrained ODE:
\begin{equation}\label{flowmap1}
y\in C_0,~~~\dot{y}=F_{0}(y)\coloneqq \begin{pmatrix}
    - \hat{f}_{q}\left(z,\chi\right)\\
    0\\
    \frac{1}{\varepsilon_p}\Psi(\omega)\chi
    \end{pmatrix},
\end{equation}
where $\hat{f}_{q}:M\times \mathbb{T}^{n}\to TM$ was defined above, and $\Psi:\R^n\to \R^{2n\times 2n}$ is the standard mapping inducing $n$-linear oscillators:
\begin{align*}
    %
    \Psi(\omega)\coloneqq \begin{pmatrix}\Omega(\omega_1)& 0 &\ldots&0\\ 0& \Omega(\omega_2)& \ldots &0\\ \vdots & \vdots & \ddots & \vdots\\ 0 & 0 & \ldots & \Omega(\omega_n)\end{pmatrix},~\Omega(\alpha) \coloneqq \begin{pmatrix}0&\alpha\\-\alpha &0 \end{pmatrix},
    \end{align*}
where $\alpha>0$. In \eqref{flowmap1}, $\hat{\chi}$ corresponds to the vector that stacks the odd components of $\chi$, $\omega_i=\tilde{\omega}_i\hat{\omega}$, where $\hat{\omega}\in \R_{>0}$ and $\tilde{\omega}_i$ is a positive rational number, and the parameters $\varepsilon_p\in\R_{>0}$ and $\varepsilon_a\in\R_{>0}$ are tunable gains. For every $q\in\mathcal{Q}$, the vector field $\hat{f}_q(z,\chi)$ is obtained by dithering the corresponding warped cost function $\tilde{\phi}_q$, formally defined below, through trajectories generated with the restricted exponential map $\text{exp}_z$. Accordingly, the dither is performed along a geodesic $\gamma:[0,1]\to M$, emanating from $z$ and which has initial velocity $\dot{\gamma}(0)$ depending on the dithering amplitudes captured by $\chi$.

To capture the switching between the different vector fields, the discrete-time dynamics $G_{0}$ of $\mathcal{H}_{0}$ are modeled by the constrained difference inclusion
\begin{align*}
    y\in D_0,~~~y^+\in G_{0}(y) \coloneqq\begin{pmatrix}
    z\\
    h(z)\\
    \chi\\
    \end{pmatrix}\tageq{\label{sES:maps}},
\end{align*}%
where the function $h:M\rightrightarrows \mathcal{Q}$, is given by
\begin{equation}
h(z) \coloneqq \big\{q\in \mathcal{Q}~:~\tilde{\phi}_q(z)=m(z)\big\},\tageq{\label{sGD:h}}
\end{equation}
and $m: M\to\R$ is defined as:
\begin{equation}\label{sGD:defm}
    m(z)\coloneqq \min_{q\in \mathcal{Q}} \tilde{\phi}_q(z).
\end{equation}
Namely, $m(z)$ is the minimum value among all the warped cost functions $\tilde{\phi}_q$ at a given point $z$. To compute $m(z)$ one only needs measurements or evaluations of $\tilde{\phi}_q(z)$, which preserves the zeroth-order nature of the hybrid algorithms. Moreover, the minimum in \eqref{sGD:defm} is well-defined since $\mathcal{Q}$ is finite, and obtaining the value of $m$ is not computationally expensive, since the complexity scales linearly with the cardinality of $\mathcal{Q}$. 

The final elements needed for the full characterization of the hybrid system $\mathcal{H}_{0}$ are the flow and jump sets $C_{0}$ and $D_{0}$, respectively. To define these sets, and since the warping induced by the diffeomorphisms is only useful if it modifies the points that are not in a neighborhood of the minimizers, we will use a threshold parameter $\gamma\in\R$, whose existence is assumed throughout the paper.

\begin{stdAssumption}\label{assumption:gammaValue}
There exists a known threshold number $\gamma \in (\underline{\phi},\phi_2)$. \QEDB
\end{stdAssumption}

Using Assumption \ref{assumption:gammaValue}, we are now ready to introduce the concept of a \emph{synergistic family of diffeomorphisms} for the solution of problem \eqref{statement:optimization}.

\begin{definition}\label{def:synergistic:family}
Let $M$ be a smooth manifold, and $\phi\in C^\infty(M)$ be a cost function satisfying Assumption \ref{assumption:phi:critical}. A family of functions $\mathcal{S}=\set{S_q}_{q\in \mathcal{Q}}$ is said to be a \emph{$\delta$-gap synergistic family of diffemorphisms adapted to $\phi$} if it satisfies:
\begin{enumerate}[label=(A$_\arabic*$)]
\item For every $q\in\mathcal{Q}$, $S_q:M\to M$ is a diffeomorphism.
\item For every $q\in\mathcal{Q}$,  $\phi(z)< \gamma \implies S_q(z) {=}z$.
\item There exists $\delta\in (0,\mu(\mathcal{S}))$, where
\begin{equation*}
    \mu\left(\mathcal{S}\right) \coloneqq \min_{\substack{p\in\mathcal{Q}\\ z\in \text{Crit}~\tilde{\phi}_q\setminus \mathcal{A}}}\max_{p\in\mathcal{Q}} \left(\tilde{\phi}_q(z) - \tilde{\phi}_p(z)\right),
\end{equation*}
and the \emph{warped cost function $\tilde{\phi}_q:M\to \R$} is defined by
\begin{equation}
    \tilde{\phi}_q \coloneqq \phi\circ S_q
\end{equation}
for all $q\in\mathcal{Q}$. \QEDB
\end{enumerate}
\end{definition}

The existence of a family of functions $\mathcal{S}$ satisfying the above properties will guarantee sufficiently enough ways of distorting the manifold $(M,g)$, so that critical points other than the minimizers of $\phi$ can be distinguished from each other. For each distortion of the manifold, a warped cost function $\tilde{\phi}_q$ can be defined, leading to a family of $N$ different vector fields in \eqref{flowmap1}. Using Definition \ref{def:synergistic:family}, we state our last main standing assumption

\begin{stdAssumption}\label{assumption:synergisticDiffeomorphisms}
There exists a $\delta$-gap synergistic family of diffeomorphisms adapted to $\phi$ with finite index set $\mathcal{Q}$. \QEDB
\end{stdAssumption}

\begin{remark}\normalfont
The verification of conditions (A1)-(A3) is clearly application-dependent, and different manifolds would usually lead to different warp costs. However, we stress that the constructions needed to implement the hybrid dynamics do not require explicit mathematical knowledge of the cost function $\phi$, but only knowledge of qualitative properties that could be verified a priori via simple tests or experiments. Particular examples of pairs $(\phi,(M,g))$ that satisfy Standing Assumptions \ref{assumption:phi:critical}-\ref{assumption:synergisticDiffeomorphisms} will be presented in Section \ref{applicartions1}. \QEDB
\end{remark}

The flow sets and jump sets of the zeroth-order hybrid dynamics can now be introduced:
\begin{align*}
    C_{0}&\coloneqq \big\{(z, q, \chi)\in M\times \mathcal{Q}\times \mathbb{T}^n
    :\left(\tilde{\phi}_q - m\right)(z)\le \delta \big\},\\
    D_{0}&\coloneqq \{(z, q, \chi)\in M\times \mathcal{Q}\times \mathbb{T}^n
    :\left(\tilde{\phi}_q - m\right)(z)\ge \delta \}.
\end{align*}
The structure of the sets $(C_{0},D_{0})$ describes the switching behavior of the hybrid dynamics $\mathcal{H}_{0}$. In particular, switches of $q$ (i.e., jumps) are allowed whenever the difference $\tilde{\phi}_q(z)-m(z)$ exceeds a $\delta$-threshold. Flows following the vector field \eqref{flowmap1} are allowed when this difference is less or equal than $\delta$. On the other hand, when the difference is exactly equal to $\delta$, flows and jumps are both allowed. This immediately indicates that solutions of $\mathcal{H}_{0}$ are not unique. Yet, the structure of the warped cost functions $\tilde{\phi}_q$ and the jump map  will prevent the existence of consecutive infinitely many jumps by inducing a hysteresis-like behavior such that whenever a solution is near a critical point of $\tilde{\phi}_q$ that is not in the set of minimizers $\mathcal{A}$, the solution will switch to a vector field that is generated from a warped cost function $\tilde{\phi}_p$ with a lower value. The existence of such a warped cost function is guaranteed by the following Lemma.
\begin{figure*}[t!]
    \centering
    \includegraphics[width=0.9\textwidth]{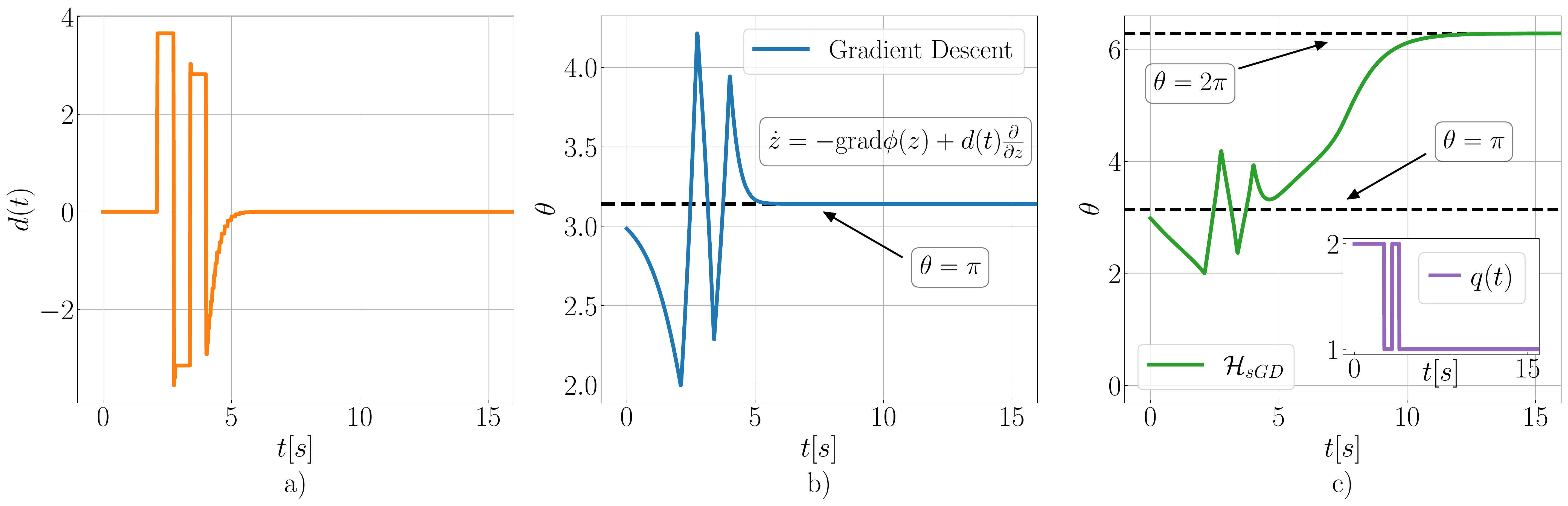}
    \caption{a) Evolution in time of additive disturbance (preserving invariance of $M$) that is able to stabilize undesired critical points in gradient flows. b) Evolution in time of trajectory of gradient flows with disturbance $d(t)\frac{\partial}{\partial z}$. c) Evolution in time $t$ of main state of $\mathcal{H}_{1}$ under the same perturbation applied to the $z$-component of the flow $F_{1}$. }
    \label{fig:disturbance:gdvshybrid}
\end{figure*}

\begin{lemma}\label{lemma:pairwise:synergistic:characterization}
Suppose that the cost function $\phi$ satisfies Assumption \ref{assumption:phi:critical}, and let $\mathcal{S}=\set{\mathcal{S}_q}_{q\in \mathcal{Q}}$ be a family of functions satisfying (A$_1$) and (A$_2$) in Definition \ref{def:synergistic:family}. If in addition $\mathcal{S}$ satisfies (A$_3$), making it a $\delta$-gap family of diffeomorphisms adapted to $\phi$, then, for every $q\in \mathcal{Q}$ and every $z\in\text{Crit}~\tilde{\phi}_q\setminus\mathcal{A}$, there exists $p\in\mathcal{Q}$ such that:
\begin{align*}
 \tilde{\phi}_p(z) + \delta < \tilde{\phi}_q(z).\tageq{\label{lemma:pairwise:synergistic}}
\end{align*}
On the other hand, if for every $q\in \mathcal{Q}$ and every $z\in\text{Crit}~\tilde{\phi}_q\setminus\mathcal{A}$, there exists $p\in\mathcal{Q}$ such that \eqref{lemma:pairwise:synergistic} is satisfied, then $\mathcal{S}$ is a $\delta$-gap synergistic family of diffeomorphisms adapted to $\phi$. \QEDB
\end{lemma}
By leveraging our standing assumptions, we now state the first main result of the paper.

\begin{theorem}\label{maintheorem}
Consider the hybrid zeroth-order dynamics $\mathcal{H}_{0}$ with data \eqref{hybridzeroorder}, and let the frequencies $\omega_i$ in \eqref{flowmap1} satisfy:
\begin{equation*}
    \tilde{\omega}_i\neq \tilde{\omega}_j, \tilde{\omega}_i\neq 2\tilde{\omega}_j, \tilde{\omega}_i\neq 3\tilde{\omega}_j,\text{ for all }i\neq j.
\end{equation*}
Then, as $(\varepsilon_p,\varepsilon_a)\to 0^+$ the set $\mathcal{A}\times \mathcal{Q}\times \mathbb{T}^n$ is GP-AS and $M\times\mathcal{Q}\times\mathbb{T}^n$ is strongly forward invariant. \QEDB
\end{theorem}

The result of Theorem \ref{maintheorem} establishes global convergence to a neighborhood of the set of minimizers $\mathcal{A}$, which can be made arbitrarily small by letting $\varepsilon_p$ and $\varepsilon_a$ go to $0^+$. Figure \ref{fig:cartoon} displays one trajectory of the hybrid system $\mathcal{H}_{0}$ that illustrates the results of Theorem \ref{maintheorem}; under sufficiently small values of $\varepsilon_p$ and $\varepsilon_a$, the trajectories of the system are constrained to evolve in the manifold $M$, and ultimately converge to the set of minimizers $\mathcal{A}$. To our best knowledge, the hybrid dynamics $\mathcal{H}_{0}$ are the first in the literature that implement deterministic zeroth-order optimization techniques with theoretical global attractivity certificates in smooth boundaryless compact Riemannian manifolds.
%
%
%
\subsection{Approximation via 1st-Order Hybrid Dynamics}
The proof of Theorem \ref{maintheorem} relies on showing that, as $(\varepsilon_{p},\varepsilon_a)\to0^+$, the trajectories of the system $\mathcal{H}_{0}$ will graphically converge to a solution of a first-order hybrid algorithm $\mathcal{H}_1$, with state $x=(z,q)$, flows given by
\begin{align}\label{sGD:flowMap}
    x\in C_1,~~~~\dot{x}&=F_{1}(x) \coloneqq \begin{pmatrix}- \nabla_{\diffp{}{z_i}}\tilde{\phi}_q(z)\diffp{}{z_i}\\
    0\end{pmatrix},
\end{align}
jumps given by the constrained difference inclusion
\begin{align}\label{sGD:jumpMap}
    x\in D_1,~~~~x^+\in G_{1}(x) \coloneqq\begin{pmatrix}
    z\\
    h(z)
    \end{pmatrix},
\end{align}
and flow set and jump set given by
\begin{subequations}\label{sGD:sets}
\begin{align}
    C_{1}\coloneqq \{(z, q)\in M\times \mathcal{Q}~:~(\tilde{\phi}_q  - m)(z)\le \delta \}\tageq{\label{sGD:flowSet}}\\
    D_{1}\coloneqq \{(z, q)\in M\times \mathcal{Q}~:~(\tilde{\phi}_q - m)(z)\ge \delta \}.\tageq{\label{sGD:jumpSet}}
\end{align}
\end{subequations}
We refer to \eqref{sGD:flowMap}-\eqref{sGD:sets} as the \emph{first-order hybrid dynamics} $\mathcal{H}_{1}:=\{C_{1},F_{1},D_{1},G_{1}\}$, since it uses first-order information of the warped costs $\tilde{\phi}_q$. Note that for every $q\in\mathcal{Q}$ the $z$-component of the vector field in \eqref{sGD:flowMap} represents a scaled version of $\text{grad}~\tilde{\phi}_q$. We note that similar vector fields have been used in the literature of optimization on manifolds in \cite{taringoo2018optimization}. In particular, the vector field differs from the coordinate representation of $$\text{grad}~\tilde{\phi}_q(z) = g^{ij}(z) \nabla_{\diffp{}{z^i}}\left(\tilde{\phi}_q\right)(z)\diffp{}{z^j},$$ by omitting the value of the Riemannian metric. Although this obstructs a coordinate-free expression of the hybrid dynamics $\mathcal{H}_{1}$, as Lemma \ref{lemma:vanishingpoints} shows, it does not modify the set of critical points of the warped cost functions.

\begin{lemma}\label{lemma:vanishingpoints}
For all $q\in\mathcal{Q}$ we have that $ \nabla_{\diffp{}{z_i}}\tilde{\phi}_q(z)\diffp{}{z_i}=0 \iff \text{grad}~\tilde{\phi}_q|_z = 0$.\QEDB
\end{lemma}

\noindent The following theorem provides a first-order version of Theorem \ref{maintheorem} for the case when the vector field \eqref{sGD:flowMap} can be explicitly computed or measured in real time, and all the standing assumptions hold.
\begin{theorem}\label{thm:stability:nominal}
The first-order hybrid dynamics $\mathcal{H}_1$ render the set $\mathcal{A}\times \mathcal{Q}$  UGAS, and the set $M\times\mathcal{Q}$ is strongly forward invariant. \QEDB
\end{theorem}

The main contribution of Theorem \ref{thm:stability:nominal} is that it overcomes the topological obstructions to global optimization on smooth compact manifolds in ODES. In particular, the asymptotic stability result is \emph{global} rather than \emph{almost} global, semi-global, or local. This result, combined with the compactness of $\mathcal{A}\times \mathcal{Q}$ and the well-posedness of the dynamics will allow us to establish important robustness properties with respect to arbitrarily small, potentially adversarial, disturbances, which could also act on the hybrid zeroth-order dynamics $\mathcal{H}_0$.
\subsection{Robustness Corollaries: Stability Under Adversarial Disturbances}
\label{secrobustness}

Since, by construction, the hybrid dynamics $\mathcal{H}_0$ and $\mathcal{H}_1$ satisfy the Basic Assumptions of \cite[Ch. 6]{bookHDS}, their stability properties are not drastically affected by small additive disturbances acting on the states and data of the hybrid systems \cite[Thm. 7.20]{bookHDS}. This observation allows us to directly obtain the following robustness corollary.
\begin{corollary}\label{corollaryrobust1}
Consider the perturbed first-order hybrid dynamics
\begin{subequations}\label{perturbed_firstorder1}
\begin{align}
&x+d_1\in C_1,~~\dot{x}=F_1(x+d_2)+d_3\\
&x+d_4\in D_1,~~x^+\in G_1(x+d_5)+d_6
\end{align}
\end{subequations}
where $\{C_1,F_1,D_1,G_1\}$ is the data of $\mathcal{H}_1$, and the signals $d_j:\text{dom}(x)\to C_1\cup D_1$, for all $j\in\{1,2,4,5,6\}$, and $d_3:\text{dom}(x)\to TC_1$, are measurable functions satisfying $$\sup_{(t,j)\in\text{dom}(x)}|d_k(t,j)|\leq d^*,$$ where $d^*>0$, for all $k\in\{1,2,\ldots,6\}$. Then, system \eqref{perturbed_firstorder1} renders the set $\mathcal{A}\times\mathcal{Q}$ GP-AS as $d^*\to0^+$. 
\QEDB
\end{corollary}

\begin{figure}[t]
    \centering
    \includegraphics[width=0.6\linewidth]{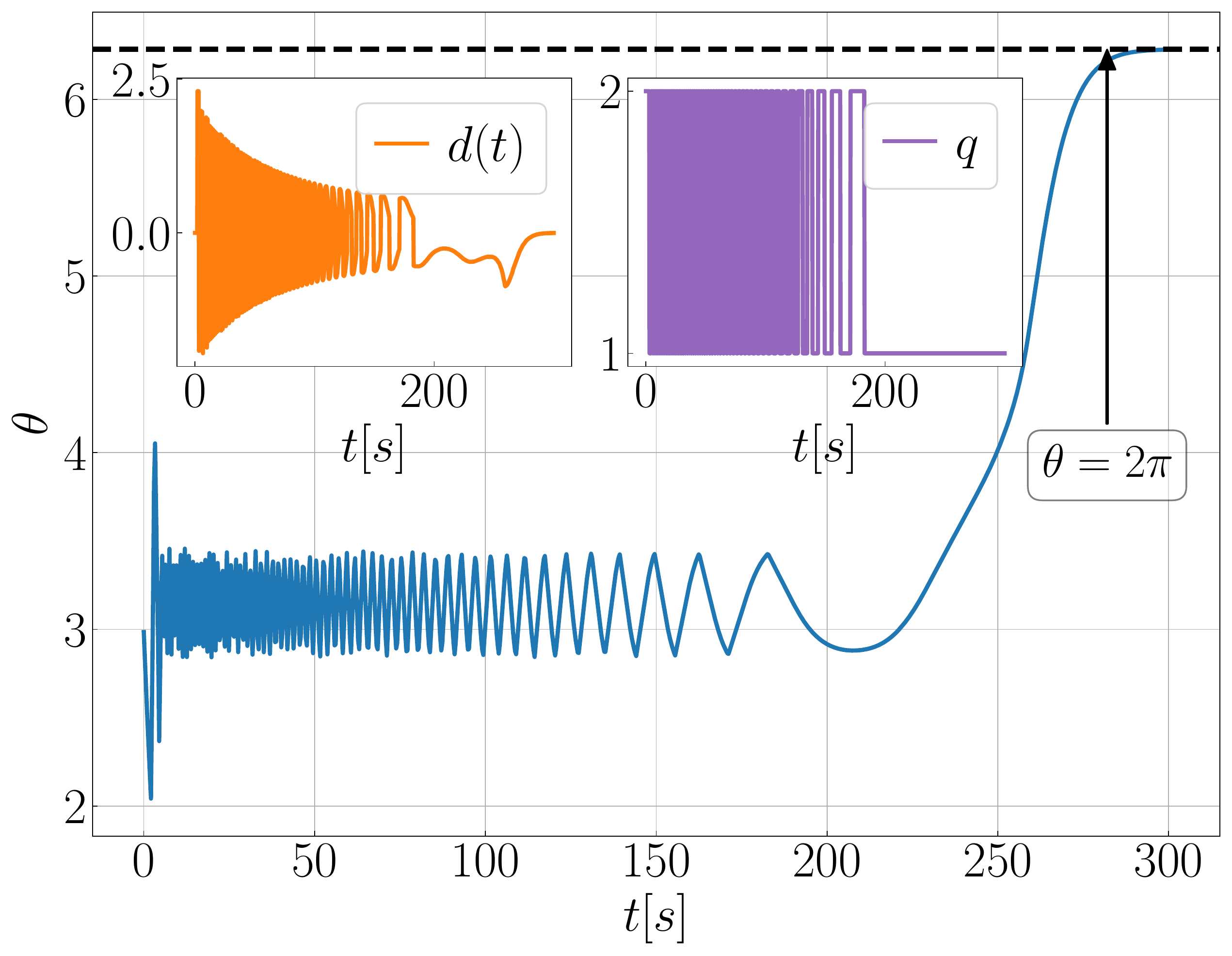}
    \caption{Evolution of trajectories of $\mathcal{H}_1$, under a small adversarial disturbance generated by a dynamical system. The insets show the amplitude of the injected disturbance, as well as the evolution of the index state $q$ in time.}
    \label{fig:DisturbanceWithHsGD}
\end{figure}
%
%

Results such as Corollary \ref{corollaryrobust1} are usually relevant for general practical applications where measurement noise or numerical approximations induce unavoidable disturbances into the system. However, the robustness properties of the hybrid controller can also be tested with respect to adversarial perturbations that are purposely designed to destabilize the set $\mathcal{A}$.  

\begin{example}\label{example36}\normalfont
Let $M=\mathbb{S}^1$ be the unit circle, which is a smooth, boundaryless compact manifold.
We consider the cost function $\phi:\mathbb{S}^1\to\mathbb{R},~z\mapsto 1-z_1$,  where $z_i\in [-1,1]$ represents the $i$-th coordinate of $z\in\mathbb{S}^1$ expressed in regular Cartesian coordinates. The cost function $\phi$ has two critical points in $\mathbb{S}^1$ corresponding to the global minimizer given by (in polar coordinates) $\theta^*=2\pi$, and a global maximizer, given by $\theta'=\pi$. We assume that the minimizer of $\phi$ is unknown, and to find this point we first implement the scaled perturbed gradient flow 
\begin{equation}\label{gradientflows}
z\in M,~~\dot{z}=- \nabla_{\diffp{}{z}}\tilde{\phi}_q(z)\diffp{}{z}+d(t)\frac{\partial}{\partial z},
\end{equation}
where $d(t)\frac{\partial}{\partial z}$ is a time-varying perturbation that preserves the invariance of $M$. This perturbation was generated by interconnecting \eqref{gradientflows} with an adversarial hybrid dynamical system to stabilize the maximizer $\theta'$. As shown in the trajectories of Figure \ref{fig:disturbance:gdvshybrid}, the adversarial perturbation is always bounded and it succeeds in stabilizing $\theta'$ (see center plot). On the other hand, when this same adversarial signal $d(t)\diffp{}{z}$ (now in open loop) is added to $\mathcal{H}_1$, as in \eqref{perturbed_firstorder1}, the hybrid dynamics guarantee convergence to the global minimizer $\theta^*$. This behavior is shown in the right plot of Figure \ref{fig:disturbance:gdvshybrid}. Finally, Figure \ref{fig:DisturbanceWithHsGD} illustrates the performance of the perturbed hybrid system under a different adversarial perturbation, which was generated using the same adversarial dynamical system used to destabilize $\theta^*$ in system \eqref{gradientflows}. In this case, the hybrid controller still guarantees convergence to $\theta^*$. \QEDB

\end{example}
An important observation from the previous discussion is that smooth model-free versions of \eqref{gradientflows} obtained via averaging might suffer from the same issues illustrated in Example \ref{example36}. In particular, if a small adversarial disturbance is able to (locally) stabilize the \emph{average dynamics} to a point not in $\mathcal{A}$, then if the disturbance is preserved by averaging, standard stability averaging results for ODEs (e.g., \cite[Ch. 10]{khalil2002nonlinear}) will predict that the actual system might also converge to a neighborhood of such a point. An example of this behavior in obstacle avoidance problems was presented in \cite[Ex. 1]{poveda2021robust}. Whether such adversarial signals can be systematically constructed in other manifolds, is an application-dependent question that we do not further pursue in this paper.

Similar to Corollary \ref{corollaryrobust1}, the following corollary provides robustness guarantees for the zeroth-order dynamics $\mathcal{H}_0$.

\begin{corollary}
Consider the perturbed zeroth-order dynamics, given by
\begin{subequations}\label{perturbed_firstorder0}
\begin{align}
&y+d_1\in C_0,~~\dot{y}_0=F_0(y+d_2)+d_3\\
&y+d_4\in D_0,~~y^+_0\in G_0(y+d_5)+d_6
\end{align}
\end{subequations}
where $\{C_0,F_0,D_0,G_0\}$ is the data of $\mathcal{H}_0$, and the signals $d_j:\text{dom}(y)\to C_0\cup D_0$, for all $j\in\{1,2,4,5,6\}$, and $d_3:\text{dom}(y)\to TC_0$, are measurable functions satisfying $$\sup_{(t,j)\in\text{dom}(y)}|d_k(t,j)|\leq d^*,$$ where $d^*>0$, for all $k\in\{1,2,\ldots,6\}$. Then, system \eqref{perturbed_firstorder0} renders the set $\mathcal{A}\times\mathcal{Q}\times\mathbb{T}^n$ GP-AS as $(d^*,\varepsilon_2,\varepsilon_1)\to0^+$. 
\QEDB
\end{corollary}
\begin{remark}
The class of problems for which smooth optimization dynamics cannot achieve robust global certificates on compact boundaryless manifolds extends beyond the case where the cost function has a unique minimizer. Indeed, by \cite[Cor. 21]{mayhew2011topological} and given a compact manifold $M$, the basin of attraction of the set of minimizers $\mathcal{A}\subset M$ of a continuous cost function $\phi$, under outer semicontinuous and locally bounded optimization dynamics $F$, is diffeomorphic to an open tubular neighborhood of $\mathcal{A}$, which, in general, is not topologically compatible with $M$. For instance, when the cost function has a finite set of global isolated minimizers $\mathcal{A} = \bigcup_{i\in I}\set{\underbar{x}_i}$, the basin of attraction is given by $\mathcal{B}_F\left(\mathcal{A}\right) = \set{x\in M~:~d(x,\underbar{x}_i) < \frac{1}{2}\min_{i\neq j}d_g\left(\underbar{x}_j,\underbar{x}_i\right)}$, which is not contractible. However, the results of Theorems \ref{maintheorem} and \ref{thm:stability:nominal} can be directly extended to overcome this type of topological obstructions for global optimization. We omit this extension due to space limitations. 
\end{remark}
\subsection{Applications: Synthesis of Algorithms}
\label{applicartions1}
We now demonstrate the effectiveness of the proposed zeroth-order hybrid dynamics $\mathcal{H}_0$ for the solution of problems of the form \eqref{statement:optimization} on two different compact Riemannian manifolds $M$. We show how to synthesize specific algorithms by generating a $\delta$-gap family of diffeomorphisms adapted to smooth cost functions defined in the unitary circle $M=\mathbb{S}^1$, and in the two-dimensional sphere $M=\mathbb{S}^2$, and we use the hybrid algorithms to achieve global model-free (practical) optimization while preserving the forward invariance of $M$. 
\subsubsection{Model-Free Feedback Optimization on $\mathbb{S}^1$}\label{sec:numExp:circle}
\begin{figure}[t!]
    \centering
    \includegraphics[width=0.65\linewidth]{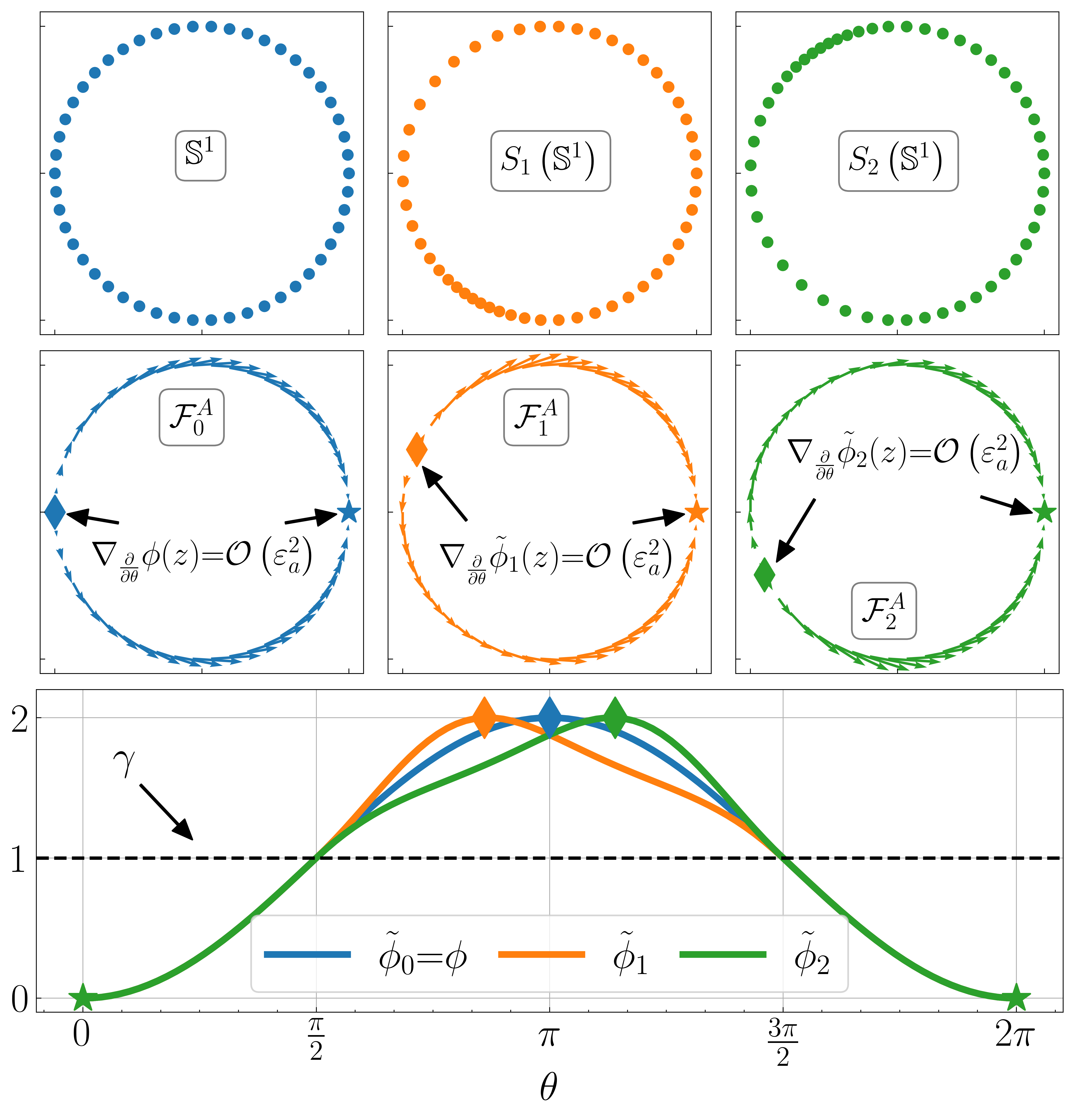}
    \caption{ (Top) Visualization of diffeomorphisms on the circle. (Middle) Average gradient-based vector fields derived from warped cost functions. (Bottom) Original and warped cost functions obtained by precomposing with diffeomorphisms. }
    \label{fig:circleDiffeomorphisms}
\end{figure}
We consider the unitary circle $\mathbb{S}^1$, and use an approach inspired by the hybrid controllers of \cite{strizic2017hybrid} to define a synergistic family of diffeomorphisms. Indeed, given  $k_q\in\R$, for $q$ belonging to some index set $\mathcal{Q}$, we define the map $S^{(1)}_q:\mathbb{S}^1\to\mathbb{S}^1$
\begin{subequations}\label{numex:circle:diffdef}
\begin{align}
    S^{(1)}_q(z)&\coloneqq\begin{cases}
             z,  &\phi(z) \le \gamma\\
             e^{k_q\alpha\left(\phi(z)-\gamma\right)\Psi}z,\quad &\phi(z)> \gamma.
            \end{cases}\\
    \Psi &\coloneqq\begin{pmatrix}0&-1\\1&0\end{pmatrix},
\end{align}
\end{subequations}
where $\alpha:\R\to\R$ is a continuously differentiable function satisfying: (B$_1$) $\alpha(0)=0$; (B$_2$) $\alpha'(0)=0$; (B$_3$) $\alpha'(r)>-1,~\forall r\ge 0$.
%
%
%
The conditions (B$_1$) and (B$_2$) ensure that $S_{q}^{(1)}$ is a continuously differentiable function that constitutes a suitable candidate for a diffeomorphism. Following the ideas of \cite{strizic2017hybrid}, we also impose condition (B$_3$) on $\alpha$. In particular, by \cite[Thm 4.1]{strizic2017hybrid}, if
\begin{figure}[t!]
    \centering
    \includegraphics[width=0.65\linewidth]{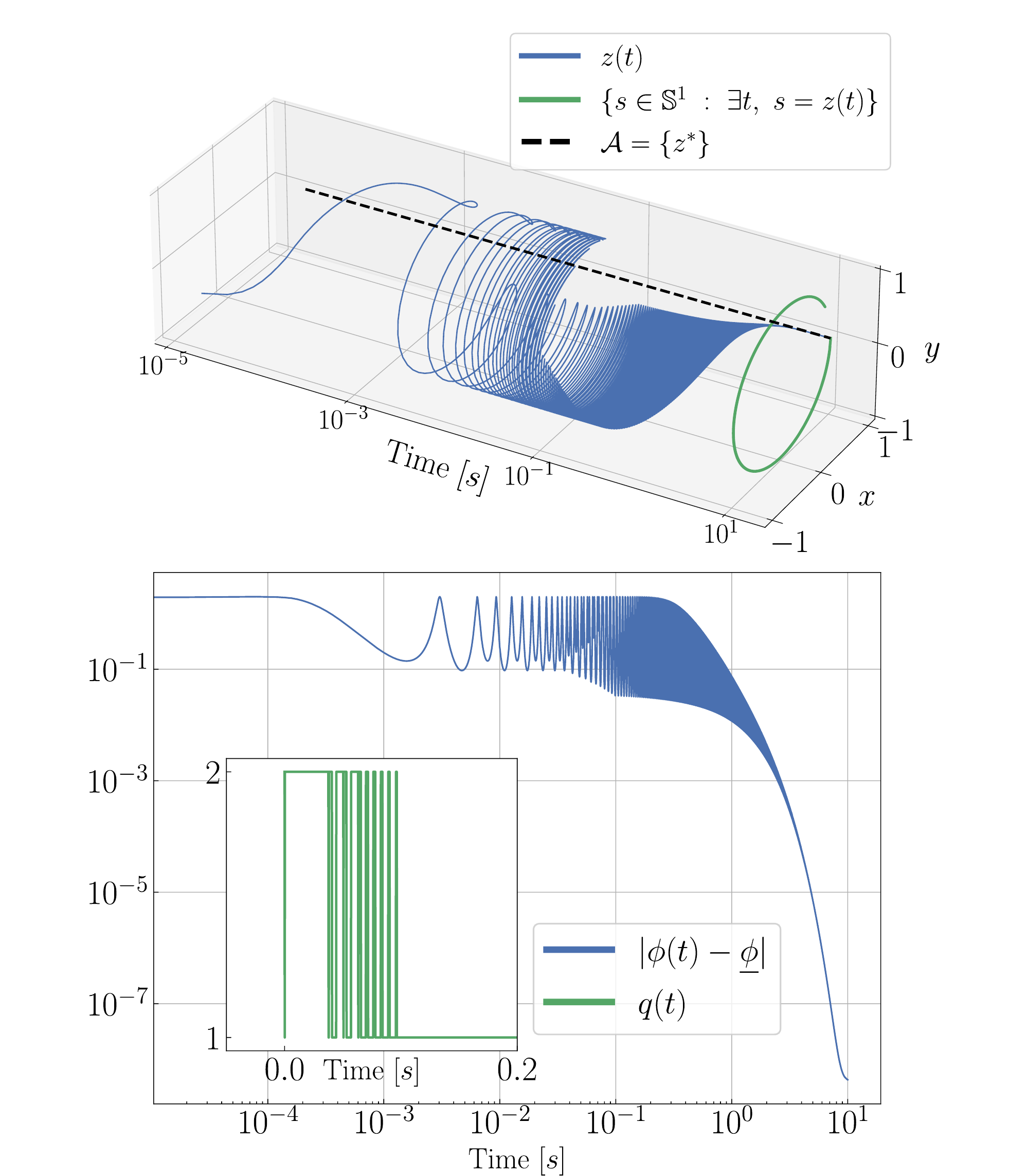}
    \caption{Model-free global optimization via $\mathcal{H}_0$ on $\mathbb{S}^1$ via Geodesic Dithering.}
    \label{fig:circle}
\end{figure}
\begin{equation*}
    \abs{k_q} < \frac{1}{\max\set{\abs{\alpha'\left(\phi(z)-\gamma\right)d\phi_z(\Psi z)}~:~z\in\mathbb{S}^1, \phi(z)\ge \gamma}},
\end{equation*}
then $S^{(1)}_q$, as defined in \eqref{numex:circle:diffdef}, is a diffeomorphism. Although the value of the bound on the admissible gains $k_q$ might not be known (since we do not know the cost function nor its differential) its existence is guaranteed by the continuity of $\alpha', \phi,$ and $d\phi$, and the compactness of $\set{z\in \mathbb{S}^1,\phi(z)\ge \gamma}$. Estimates of the bound could be obtained by a Monte Carlo method that uses measurements of $\phi$ and $d\phi$ at different $z\in\mathbb{S}^1$.

Given a cost function $\phi^{(1)}: \mathbb{S}^1\to \R$, and choosing an appropriate set of gains $\set{k_q}_{q\in\mathcal{Q}}$ with corresponding diffeomorphisms defined by \eqref{numex:circle:diffdef}, it is possible to build a suitable $\delta$-gap synergistic family of diffeomorphisms subordinate to $\phi^{(1)}$. To illustrate this process, and as in Example \ref{example36}, consider again the cost function defined by $\phi^{(1)}(z)\coloneqq 1- z_1$. Assume that only measurements from $\phi^{(1)}$ are available for the purpose of feedback design, but that the intermediate value $\gamma^{(1)} \coloneqq 1 \in (0,2)= \left(\underline{\phi}^{(1)}, \phi^{(1)}_2\right)$ and the number of critical points of $\phi^{(1)}$ are known in advance. Let $\alpha(r)=r^2$, and note that it satisfies conditions (B$_1$)-(B$_3$). Then, choosing  any two (number of critical points of $\phi$) different gains satisfying the bound on $|k_q|$, we can obtain a synergistic family of diffeomorphisms subordinate to $\phi^{(1)}$. Indeed, choosing $\mathcal{Q}=\set{1,2},\abs{k_q}<1,~q\in \mathcal{Q},k_1\neq k_2$, the set $\mathcal{S}^{(1)}=\set{S^{(1)}_q}_{q\in\mathcal{Q}}$ is a $\delta$-gap family of diffeomorphisms adapted to $\phi^{(1)}$ with gap $\delta <\mu\left(\mathcal{S}^{(1)}\right)$. In Figure \ref{fig:circleDiffeomorphisms} we show a visualization of the diffeomorphisms in this family using the choice $k_1=\frac{1}{2},k_2=-\frac{1}{2}$. We show how these diffeomorphisms warp the original cost function, and we also plot the gradient-based vector fields obtained from the warped cost functions, corresponding to the average-dynamics of the zeroth-order algorithm $\mathcal{H}_{0}$. 

By implementing $\mathcal{H}_{0}$ on the unitary circle, we obtain the results shown in Figure \ref{fig:circle}. As shown in the plots, the zeroth-order hybrid dynamics with geodesic dithering are able to converge to the minimizer of $\phi^{(1)}$, i.e. to $\underline{z} = (1,0)$ while avoiding the other critical point $\overline{z} = (-1,0)$, and maintaining invariance of the manifold. The dynamics are implemented using zeroth-order information from the cost function $\phi^{(1)}$.
\subsubsection{Model-free Feedback Optimization on $\mathbb{S}^2$}
We now consider the two-dimensional sphere $\mathbb{S}^3$. Following similar ideas as in Section \ref{sec:numExp:circle}, given  $k_q\in\R$, for $q\in\mathcal{Q}$, and $u\in \mathbb{S}^2$, we define the map $S^{(2)}_{q,u}:\mathbb{S}^2\to\mathbb{S}^2$ as follows:
\begin{subequations}\label{numex:sphere:diffdef}
\begin{align}
&S^{(2)}_{q,u}\coloneqq \begin{cases}
             z,  &\phi(z) \le \gamma,\\
            e^{k_q\alpha\left(\phi(z)-\gamma\right)[u]_\times}z,~ &\phi(z)> \gamma.
            \end{cases}\\
[u]_\times &: \R^3 \to \R^{3\times 3}\\
&u\longmapsto \begin{pmatrix}0&-u_3&u_2\\u_3&0&-u_1\\-u_2&u_1&0\end{pmatrix},
\end{align}
\end{subequations}
where $\alpha:\R\to\R$ is again a continuously differentiable function satisfying (B$_1$)-(B$_3$). The definition of the map results from modifiying the function introduced in \cite[Sec 3.4.3]{mayhew2011synergistic}, for the angular warping of the two-dimensional sphere, by using the function $\alpha$ and letting the warping act only when $\phi$ exceeds the known threshold $\gamma$. Continuous differentiability of the resulting map follows directly by the conditions imposed on $\alpha$. The following lemma extends \cite[Thm. 4.1]{strizic2017hybrid} to the 2-sphere:

\begin{lemma}\label{lemma:kq:s2}
Let $k_q$ satisfy the following inequality:
\begin{equation}\label{s2:lemma:kq}
    \abs{k_q} {<} \frac{1}{\max\set{\abs{\alpha'\left(\phi(z){-}\gamma\right)d\phi_z\left([u]_\times z\right)}:z\in\mathbb{S}^2, \phi(z)\ge\gamma}},
\end{equation}
then $S^{(2)}_{q,u}$, as defined in \eqref{numex:sphere:diffdef}, is a diffeomorphism.\QEDB
\end{lemma}
%
\begin{figure*}[t!]
    \centering
    \includegraphics[width=\linewidth]{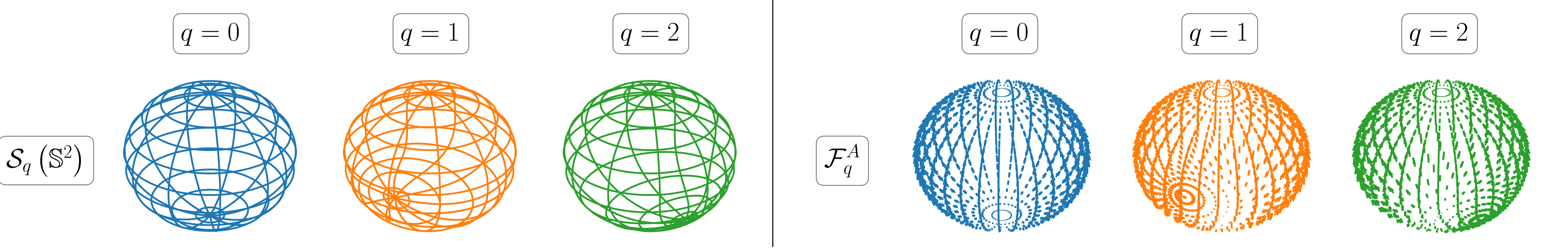}
    \caption{Left) Visualization of diffeomorphisms on the sphere. (Right) Average gradient-based vector fields derived from warped cost functions}
    \label{fig:sphereDiffeomorphisms}
\end{figure*}
\textbf{Proof:} The proof follows a similar structure to the one presented in \cite[App. 7]{mayhew2010hybrid}. First, we show that the map is proper. Indeed, recall that for every well defined Riemannian metric $g$ on $\mathbb{S}^2$, it follows that $\mathbb{S}^2$ is a metric space with the metric given by the distance function $d_g$ defined in \eqref{riemannian:distance}. Since every metric space is Hausdorff, and $\mathbb{S}^2$ is compact, by its continuity, the map $S_{q,u}^{2}$ is proper via \cite[Prop. 4.93 a]{lee2010introduction}.\\
On the other hand, we compute the differential of the map. When $\phi(z)\leq \gamma$, it follows trivially that $d\left(S^{(2)}_{q,u}\right)_z=I$. Now, whenever $\phi(z)> \gamma$, by using the fact that $$d\left(e^{A\eta(z)}\right)_z= e^{A\eta(z)}\left(I+Azd\eta_z\right)$$ for $A\in\R^{n\times n}$ and $\eta:\R^n\to \R$, the differentiability of $\alpha:\R\to \R$, and the chain rule of differentiation, the differential of $S^{(2)}_{q,u}$ satisfies:
\begin{align}
    d\left(S^{(2)}_{q,u}\right)_z &=  e^{k_q\alpha\left(\phi(z)-\gamma\right)[u]_\times}\big(I \notag\\
    &\quad\qquad+ k_q\alpha'\left(\phi(z)-\gamma\right)[u]_\times z d\phi_z\big).\tageq{\label{sec:num:jacobianMatrix}}
\end{align}
By using the fact that $e^{k_q\alpha\left(\phi(z)-\gamma\right)[u]_\times}$ is a unitary rotation satisfying $\det\left(e^{k_q\alpha\left(\phi(z)-\gamma\right)[u]_\times}\right)=1$, it follows that the Jacobian determinant is given by
\begin{equation*}
\det\left(d\left(S^{(2)}_{q,u}\right)_z \right) = 1+ k_q\alpha'\left(\phi(z)-\gamma\right)d\phi_z\left([u]_\times z\right),
\end{equation*}
which never vanishes if \eqref{s2:lemma:kq} is satisfied. Therefore, since $\mathbb{S}^2$ is simply connected,  the map is proper and the Jacobian determinant of the map never vanishes, it follows that $S^{(2)}_{q,u}$ is a diffeomorphism via \cite[Thm. B]{gordon1972diffeomorphisms}. \hfill $\blacksquare$

\begin{figure}[t]
    \centering
    \includegraphics[width=0.45\linewidth]{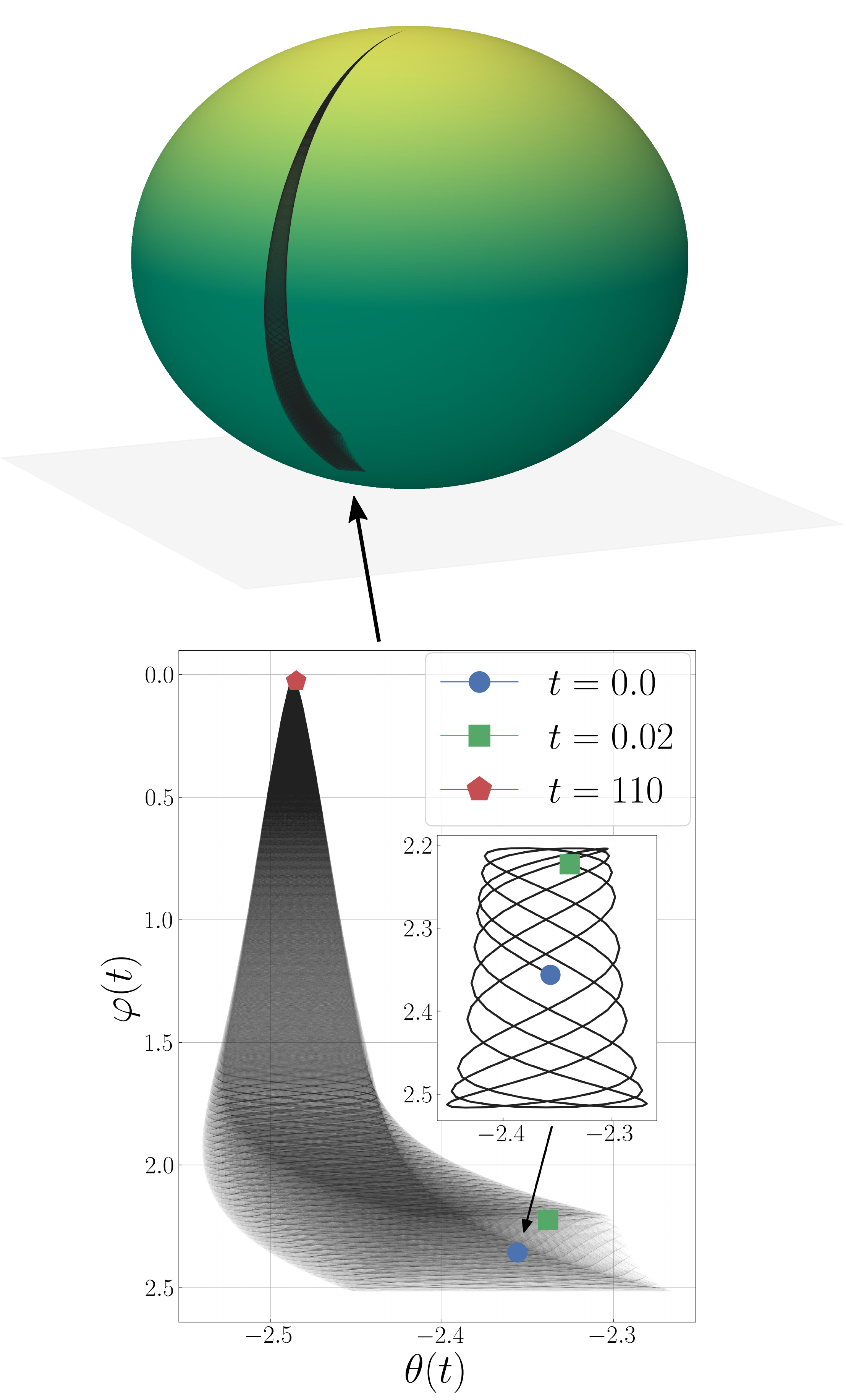}
    \caption{Synergistic Model-Free Optimization Seeking on $\mathbb{S}^2$ via Geodesic Dithering.}
    \label{fig:sphere}
\end{figure}

Now, consider the cost function $\phi^{(2)}(z) = 1-z_3$ evolving in $\mathbb{S}^2$, and where $z_i\in [-1,1]$ represents the $i$-th coordinate of $z\in\mathbb{S}^2$ expressed in regular Cartesian coordinates. By assuming that the threshold value $\gamma = 1\in \left(\underline{\phi}^{(2)}, \phi_2^{(2)}\right)$ is previously known, and that the number off critical points of $\phi^{(2)}$ is two, we can pick the gains $k_1=\frac{1}{2},k_2=-\frac{1}{2}$, the vector $u= (0,1,0)\in \R^3$, and the associated family of maps $\mathcal{S}^{(2)}\coloneqq\set{S^{(2)}_{q,u}}_{q\in\mathcal{Q}}$ with $\mathcal{Q}\coloneqq\set{1,2}$ as a candidate family of functions for warping the space. Since $\abs{k_q}<1$ for all $q\in\mathcal{Q}$, via Lemma \ref{lemma:kq:s2}, $\mathcal{S}^{(2)}$ is a family of diffeomorphisms. In fact, by Lemma \ref{lemma:pairwise:synergistic:characterization}, $S^{(2)}$ is a $\delta$-synergistic family of diffeomorphisms adapted to $\phi$ with gap $\delta<\frac{1}{4}$. Figure \ref{fig:sphereDiffeomorphisms} shows a visualization of the diffeomorphisms in this family with the choice $k_1=\frac{1}{2},k_2=-\frac{1}{2}$. The figure also shows the vector fields obtained from the warped cost functions.\\
Using $\mathcal{S}^{(2)}$, we implement the HDS $\mathcal{H}_{0}$ and obtain the results shown in Figure \ref{fig:sphere}. In the figure, we zoom in on the solution, and we plot the trajectory of the state in the azymuth vs. elevation angles phase space, showing how the dithering is compliant with the geodesics on the sphere. The trajectory shows that the state $z$ converges to the global minimizer $z^* = (0,0,1)$, while avoiding the critical point $\overline{z}=(0,0,-1)$.
\section{Analysis and Proofs of Main Results}
\label{sec_analysis}
In this section, we present the proofs of our main results. Since the stability results of the zeroth-order hybrid dynamics $\mathcal{H}_0$ in Theorem \ref{maintheorem} rely on the stability properties of the first-order dynamics $\mathcal{H}_1$, we first present the proof of Theorem \ref{thm:stability:nominal}.
\subsection{Proof of Theorem \ref{thm:stability:nominal}}
We first present the proof of our auxiliary lemmas.

\textbf{Proof of Lemma \ref{lemma:pairwise:synergistic:characterization}:}
    Suppose that $\mathcal{S}$ is a $\delta$-gap synergistic family of diffeomorphisms adapted to $\phi$. Then, we have that $\delta<\mu(\mathcal{S})$, meaning that
    \begin{align}
        &\delta < \max_{p\in\mathcal{Q}} \left(\tilde{\phi}_q - \tilde{\phi}_p\right)(z)\quad \forall q\in\mathcal{Q},\tageq{\label{lemma:pairwise:synergistic:aux}},
    \end{align}
    for all $z\in\text{Crit}~\tilde{\phi}_q\setminus \mathcal{A}.$
    From \eqref{lemma:pairwise:synergistic:aux} it follows directly that for all $q\in \mathcal{Q}$ and $z\in\text{Crit}~\tilde{\phi}_q\setminus \mathcal{A}$, there exists $p\in \mathcal{Q}$ such that \eqref{lemma:pairwise:synergistic} is satisfied.\\
    Conversely, assume that for every $q\in \mathcal{Q}$ and $z\in\text{Crit}\left(\tilde{\phi}_q\right)\setminus\mathcal{A}$, there exists $p\in\mathcal{Q}$ such that \eqref{lemma:pairwise:synergistic} is satisfied. In particular, for all $q\in\mathcal{Q}$ and $z\in\text{Crit}\left(\tilde{\phi}_q\right)\setminus\mathcal{A}$ it follows that
    $
        \delta <  \max_{p\in\mathcal{Q}} (\tilde{\phi}_q - \tilde{\phi}_p)(z),
    $
    which in turn implies that
    \begin{align*}
        \delta <  \min_{\substack{q\in\mathcal{Q}\\ z\in \text{Crit}~\tilde{\phi}_q\setminus \mathcal{A}}}\max_{p\in\mathcal{Q}} (\tilde{\phi}_q - \tilde{\phi}_p)(z).
    \end{align*}
    This concludes the proof. \hfill $\blacksquare$

\begin{lemma}\label{lemma:wellPosed}
The HDS $\mathcal{H}_{1}$ is a well-posed hybrid dynamical system.
\end{lemma}

\textbf{Proof:} We prove that $\mathcal{H}_{1}$ satisfies the hybrid-basic conditions \cite[Def. 2.20]{sanfelice2020hybrid}. First, note that $F_{1}:M\times \mathcal{Q}\to T_zM$ is continuous. Indeed, let $U\subset T_zM$ be open, and consider $$\tilde{f}^{-1}(U) = \set{(z,q)\in M\times\mathcal{Q}~:~\nabla_{\diffp{}{z_i}}\tilde{\phi}_q(z)\diffp{}{z_j}\in U}.$$ Because $\mathcal{Q}$ is a finite discrete set endowed with the discrete topology, and due to the fact that open sets in $M\times \mathcal{Q}$ are sets of the form $U_M\times U_{\mathcal{Q}}$, with $U_M$ and $U_{\mathcal{Q}}$ open in $M$ and $\mathcal{Q}$ respectively, in order to prove that $\tilde{f}^{-1}(U)$ is open, and hence that $\tilde{f}$ is continuous, it suffices to show that $\pi_{M}\left(\tilde{f}^{-1}\left(U\right)\right)$ is open. To do so, define $\tilde{\mathcal{Q}}\coloneqq \set{q\in\mathcal{Q}~:~ \exists (z,u)\in M\times U\text{ with } \nabla_{\diffp{}{z_i}}\tilde{\phi}_q(z)\diffp{}{z_j} =u }$ and notice that:
\begin{align*}
    \pi_{M}\left(\tilde{f}^{-1}(U)\right) = \bigcup_{q\in\tilde{\mathcal{Q}}} \left(\nabla_{\diffp{}{z_i}}\tilde{\phi}_q(z)\diffp{}{z_i}\right)^{-1}\left(U\right)\tageq{\label{lemma:wellPosed:aux}}.
\end{align*}
Since for every fixed $q\in\mathcal{Q}$ it follows that $\nabla_{\diffp{}{z_i}}\tilde{\phi}_q(z)\diffp{}{z_i}$ is continuous by smoothness of $\phi$ and the fact that $S_q$ is a diffeomorphism, we obtain that $\left(\nabla_{\diffp{}{z_i}}\tilde{\phi}_q(z)\diffp{}{z_i}\right)^{-1}\left(U\right)$ is open for all $q\in\mathcal{Q}$. Hence, via \eqref{lemma:wellPosed:aux}, $\pi_{M}\left(\tilde{f}^{-1}(U)\right)$ is open, meaning that $-\nabla_{\diffp{}{z_i}}\tilde{\phi}\diffp{}{z_i}$ is continuous in $M\times Q$ and consequently that $F_{1}$ is also continuous. Moreover, since $M$ is compact by assumption and $\mathcal{Q}$ is trivially compact, by continuity, we have that $F_{1}$ is bounded on $M\times \mathcal{Q}$.\\
On the other hand,  define $\tilde{\mu}:M\times \mathcal{Q}\to \R$ by $\tilde{\mu}(z,q)\coloneqq (\tilde{\phi}_q - m)(z) $, and note that it is continuous by similar arguments to the ones presented to prove continuity of $F_{1}$. Moreover, since $\text{gph}~ h = \set{(z,q)\in M\times \mathcal{Q}~:~ z\in M,~\tilde{\mu}(z,q)=0}$ is closed due to the continuity of $\tilde{\mu}$, it follows that $h$, and hence, that $G_{1}$ are outer-semicontinuous. Boundedness of $G_{1}$ thus follows by compactness of $M\times Q$.\\
We conclude by noting that $C_{1}$ and $D_{1}$ are closed, since they are sublevel and superlevel sets, respectively, of the continuous function $\tilde{\mu}$. The result follows via \cite[Thm. 6.30]{bookHDS}. \hfill $\blacksquare$
\textbf{Proof of Lemma \ref{lemma:vanishingpoints}}:
    Consider an arbitrary $\phi\in C^\infty(M)$ and assume that $\grad{\phi}|_z = 0$ at some $z\in M$. Then, by the coordinate representation of $\grad{\phi}$, it follows that:
    $
         g^{ij}(z)\nabla_{\diffp{}{z_i}}\phi(z) \diffp{}{z_j} =0.
    $
    Therefore, since $g^{ij}(z)$ is nonsingular for all $z\in M$, and $\set{\diffp{}{z_i}}_{i=1}^n$ constitutes a local frame, it follows that
    \begin{align*}
        \nabla_{\diffp{}{z_i}}\phi(z) =0,~\forall i\in\set{1,\cdots, n},\tageq{\label{lemma:vanishingpoints:aux}}
    \end{align*}
    which implies that $ \nabla_{\diffp{}{z_i}}\phi(z)\diffp{}{z_i} =0$. Conversely, assume that $ \nabla_{\diffp{}{z_i}}\phi(z)\diffp{}{z_i} =0$. Then, \eqref{lemma:vanishingpoints:aux} is satisfied, and thus it follows that $0= g^{ij}(z)\nabla_{\diffp{}{z_i}}\phi(z) \diffp{}{z_j} = \grad{\phi}|_z.$ This concludes the proof. \hfill $\blacksquare$
\noindent Now, let
\begin{align*}
    \mathcal{E} \coloneqq \set{(z,q)\in M\times \mathcal{Q}~:~ z\in \text{Crit}~\tilde{\phi}_q\setminus \mathcal{A} },
\end{align*}
represent the set of all the critical points of the cost function, warped by the family of diffeomorphisms, that are not minimizers of $\phi$. The following lemma shows that $\mathcal{E}$ is properly contained in $D_{1}$, and thus, that jumps are enforced whenever $(z,q)\in \mathcal{E}$. Therefore, due to the definition of $C_{1}$, this implies that no flow is allowed in critical points that are not minimizers of $\phi$.
\begin{lemma}\label{lemma:wellDefined:jumpSet}
Suppose that Assumption \ref{assumption:synergisticDiffeomorphisms} is satisfied. Then $\mathcal{E}\subset D_{1}^\circ$ and $G_{1}(\mathcal{E})\subset C_{1}^\circ$.
\end{lemma}
\textbf{Proof:} The proof follows directly from  Lemma \ref{lemma:pairwise:synergistic:characterization}. Indeed, assume that $(z,q)\in \mathcal{E}$, and hence that $z\in \text{Crit}~\tilde{\phi}_q\setminus \mathcal{A} $. By the fact that $\mathcal{S}$ is a $\delta-$gap synergistic family of diffeomorphisms adapted to $\phi$, and using Lemma \ref{lemma:pairwise:synergistic:characterization}, there exists $p\in\mathcal{Q}$ such that $\tilde{\phi}_p(z) + \delta < (\tilde{\phi}_q)(z)$. Since by definition $m(z) \leq (\tilde{\phi}_q)(z)$ for all $p\in\mathcal{Q}$ we additionally obtain that
    \begin{align*}
       m(z) + \delta < \tilde{\phi}_q(z),
    \end{align*}
    which implies that $(z,q)\in D_{1}$, and thus that  $\mathcal{E}\subseteq D_{1}$.\\ Now, let $(z,q)\in \partial D_{1}$, i.e. be such that $m(z) = \tilde{\phi}_q(z) - \delta$. For the sake of contradiction assume  that $(z,q)\in\mathcal{E}$. Again, by Lemma \ref{lemma:pairwise:synergistic:characterization}, there would exist $p\in\mathcal{Q}$ such that $$\tilde{\phi}_p(z) < \tilde{\phi}_q(z) - \delta  = m(z).$$ But this contradicts the fact that $m(z)\leq \tilde{\phi}_q(z)$ for all $q\in\mathcal{Q}$. Hence $(z,q)$ cannot be in $\mathcal{E}$, which in turn shows that there exists an element of $\partial D_{1}$ that is not in $\mathcal{E}$, and thus that $\mathcal{E}\subset D_{1}^\circ$.\\
    The fact that $G_{1}(\mathcal{E})\subset C_{1}^\circ$ follows by construction since after a jump we have that $\tilde{\phi}_{q^+}(z^+) - m(z^+) =\tilde{\phi}_{q^+}(z) - m(z) = 0 < \delta$. This concludes the proof. \hfill $\blacksquare$

By leveraging the results of the previous lemmas we are prepared to prove the first main theorem of the paper.

\noindent 
\textbf{Proof of Theorem \ref{thm:stability:nominal}:} Consider the Lyapunov function candidate
\begin{align*}
    V(x) &\coloneqq \tilde{\phi}_q(z) - \underline{\phi},
\end{align*}
which is continuous due to similar arguments to the ones used to prove continuity of $F_{1}$ in Lemma \ref{lemma:wellPosed}.
Due to the fact that $\underline{\phi}<\phi(z)$ for all $z\not\in \mathcal{A}$, in conjunction with (A$_2$) in Definition \ref{def:synergistic:family}, we have that $\tilde{\phi}_q(z) - \underline{\phi} \ge 0$ for all $(z,q)\in M\times \mathcal{Q}$ and $\tilde{\phi}_q(z) - \underline{\phi} = 0$ if and only if $(z,q)\in\mathcal{A}\times \mathcal{Q}$. Therefore, it follows that $V(x)\ge 0$ for all $(z,q)\in M\times \mathcal{Q}$ and $V(x)=0$ if and only if $z\in\mathcal{A}$.\\
Now, during the flows of $\mathcal{H}_{1}$ the change of the Lyapunov function $\mathcal{L}_{F_{1}} V (x) = dV_x\left(F_{1}(x)\right)$ is given by:
\begin{align*}
    dV_x\left(F_{1}(x)\right)
                                &= -\diffp{\tilde{\phi}_q}{z_j}[z]dz^j\left(\nabla_{\diffp{}{z^i}}\tilde{\phi}_q(z)\diffp{}{z^i}\right),\tageq{\label{thm:stability:nominal:aux}}
\end{align*}
where we have used the linearity of the differential, and the local coordinate representation of $d\left(\tilde{\phi}_q\right)_z = \diffp{\tilde{\phi}_q}{z_i}[z]dz_i$ at $z\in M$. Using linearity of $dz_i$, the fact that $\nabla_X f = Xf$ for all $X\in \mathfrak{X}(M), f\in C^\infty(M)$, and $dz_j\left(\diffp{}{z_i}\right) = \delta^{j}_i$, from \eqref{thm:stability:nominal:aux}, we get:
\begin{align*}
    dV_x\left(F_{1}(x)\right)
    &=u_C(x), \quad \forall  x\in C_{1},\tageq{\label{thm:stability:nominal:flows}}
\end{align*}
where $u_C:M\times Q\to \R$ is defined as 
$
    u_C(x) \coloneqq -\sum_{i=1}^n\left(\diffp{\tilde{\phi}_q}{z_i}[z]\right)^2.
$
Moreover, during jumps of $\mathcal{H}_{1}$:

\begin{align*}
    V\left(x^+\right) - V(x)
                &= \left(\tilde{\phi}_{h(z)}-\tilde{\phi}_q\right)(z)\\
                &= \left(m - \tilde{\phi}_q\right)(z)\\
                &\leq u_D(x),\quad\forall x\in D_{1},\tageq{\label{thm:stability:nominal:jumps}}
\end{align*}
with $u_D:M\times \mathcal{Q}\to \R$ defined by $u_D(x)\coloneqq -\delta$, and where we have used the fact that $$(z,q)\in D_{1}\implies \left(\tilde{\phi}_q-m\right)(z) \ge \delta.$$ In particular, due to the fact that $\mathcal{A}\subset C_{1}^\circ$ and $u_{C}(z,q) = 0$ for all $(z,q)\in\mathcal{A}$, from \eqref{thm:stability:nominal:flows} we have that $\mathcal{A}$ is stable under $\mathcal{H}_{1}$.\\
Now, in order to show attractivity of $\mathcal{A}$  we employ the hybrid invariance principle \cite[Thm. 3.23]{sanfelice2020hybrid}. Indeed, since $u_C(x)\leq 0$ for all $x\in C_{1}$ and $u_D(x)< 0$ for all $x\in D_{1}$, and using $u_D^{-1}(0)=\emptyset$, given $r\in V\left(\mathcal{A}\cup \mathcal{E}\right) \subset [0,\overline{\phi}-\underline{\phi}]$, solutions approach the largest weakly invariant set in
\begin{align*}
    V^{-1}(r) \cap \left(\mathcal{A}\cup \mathcal{E}\right) .\tageq{\label{thm:nominal:invariance}}
\end{align*}
Let $\Omega$ denote such largest invariant set and first assume that $r\neq 0$. By the definition of $\mathcal{E}$ and the synergistic family of diffeomorphisms, it follows that $\Omega\subset \mathcal{E}$. Moreover, by Lemma \ref{lemma:wellDefined:jumpSet}, we obtain that $\Omega\subset D_{1}^\circ$. Since $D_{1}^\circ \cap C_{1}=\emptyset$ by construction, for $\Omega$ to be invariant under $\mathcal{H}_{1}$, we would need to have that $\Omega = G(\Omega)$, but this would imply, via Lemma \ref{lemma:wellDefined:jumpSet}, that $\Omega\subset C_{1}^\circ$, and thus that $\Omega \subset C_{1}^\circ \cap D_{1}^\circ  = \emptyset\implies \Omega = \emptyset$. Therefore, we necessarily have that $r=0$, and thus that $\forall (z(0),q(0))\in M\times \mathcal{Q}$ solutions approach the largest weakly invariant set in $V^{-1}(0)\cap \left(\mathcal{A}\cup\mathcal{E}\right) = \mathcal{A}$, which is $\mathcal{A}$ it self. UGAS follows directly by the global attractivity and stability of $\mathcal{A}$. \hfill $\blacksquare$
\subsection{Proof of Theorem \ref{maintheorem}}
 The proof uses tools recently developed for averaging on compact Riemannian manifolds \cite{taringoo2018optimization} in conjunction with the framework for hybrid extremum seeking control introduced in \cite{poveda2017framework,Poveda2021RobustHZ}.
 
In particular, since $M$ is compact, we can select $\varepsilon_a\in\R_{>0}$ such that $\text{exp}_{z}\left( \varepsilon_a\hat{\chi}^i\diffp{}{z_i}\right)\in \imath (M)$, with $\imath(M)$ the injectivity radius of $M$ \cite[Lemma 3.2]{taringoo2018optimization}. This makes possible a  Taylor expansion in normal coordinates along the geodesic dithers for every $\tilde{\phi}_q$, such that the average dynamics of  $\mathcal{H}_{0}$ are given by $\mathcal{H}_{0}^A = \left(C_1, F_{0}^A, D_{1}, G_{1}\right)$, where $C_{1},D_{1}$, and $G_{1}$ are defined in \eqref{sGD:flowSet}, \eqref{sGD:jumpSet} and \eqref{sGD:jumpMap} respectively, and $F_{0}^A:M\times \mathcal{Q}\to TM\times \mathbb{N}$ is:
\begin{align*}
    F_{0}^A(x) & {\coloneqq} \begin{pmatrix}\mathcal{F}^{A}_q(z)\\\mathfrak{F}^{A}(q)\end{pmatrix} {=} \begin{pmatrix} -\nabla_{\diffp{}{z_i}}\tilde{\phi}_q\diffp{}{z_i}~ {+} \sum_{i=1}^n\mathcal{O}_q(\varepsilon_a^2)\diffp{}{z_i}\\ 0 \end{pmatrix},
\end{align*}
where $q$ in $\mathcal{O}_q$ indicates the dependence of the residual error on each particular diffeomorphism.
Hence, it follows that, on closed subsets of $M$, we have that
\begin{equation}\label{thm:sgpas:regularization}
F_{0}^A(x)\in \overline{\text{con}}_zF_{1}(x+k\varepsilon_a^2\mathbb{B},0)+\left(k\varepsilon_a^2\mathbb{B},0\right),
\end{equation}
for some $k>0$, where $F_{1}$ was defined in \eqref{sGD:flowMap}. In \eqref{thm:sgpas:regularization} the convex hull is done on the $z$ variable only, and the Minkowski additions ($z+k\varepsilon_a^2\B$) are performed in a suitable ambient euclidean space which always exist due to the Whitney Embedding Theorem \cite[Thm 6.15]{lee2013smooth}. 
Now, due to \eqref{thm:sgpas:regularization}, any solution of the average dynamics $\mathcal{H}^A_{0}$ is also a solution of an inflated HDS generated from $\mathcal{H}_{1}$. Hence, and since $\mathcal{H}_{1}$ is a well-posed HDS via Lemma \ref{lemma:wellPosed}, by \cite[Thm. 7.21]{bookHDS} we conclude that system $\mathcal{H}^A_{0}$ renders the set $\mathcal{A}$ GP-AS in the ambient euclidean space as $\varepsilon_a\to0^+$. Since $\mathcal{H}^A_{0}$ and $\mathcal{H}_{1}$ are nominally well-posed, all conditions to apply \cite[Cor. 1]{poveda2017framework} are satisfied. Therefore, in conjunction with the compactness of $M$, $\mathcal{H}_{0}$ renders the set $\mathcal{A}\times \mathbb{T}^n$ GP-AS in the ambient euclidean space as $(\varepsilon_p,\varepsilon_a)\to0^+$. Note that any trajectory of the state $z$ involved in $\mathcal{H}_{0}$ is constrained to evolve on $M$ since the dithering is performed along geodesics on the manifold, and due to the fact that $\hat{f}_q(z,\chi)\in T_zM$ for all $(z,q,\chi)\in M\times \mathcal{Q}\times \mathbb{T}^n$. Therefore, GP-AS of $\mathcal{A}$ is in fact in the sense of Definition \ref{definitionSGPAS}. This obtains the result. \hfill $\blacksquare$
\section{Conclusions and Outlook}\label{sec:conclusions}
We introduced a class of zeroth-order hybrid algorithms for the global solution of model-free optimization problems on smooth, compact, boundaryless manifolds. The hybrid algorithms combine continuous-time dynamics and discrete-time dynamics to achieve robust \emph{global} practical stability of the optimizer of a smooth cost function accessible only via measurements. In this way, the dynamics can overcome the topological obstructions that prevent the solution of this problem using algorithms characterized by smooth ODEs. The stability and robustness properties of the algorithms were characterized using tools from hybrid dynamic inclusions. By leveraging these robustness properties, future research directions will explore slowly time-varying optimizers as well as the optimization of the steady state of plants in the loop.
A completely coordinate-free formulation of the hybrid dynamics, as well as the development of accelerated dynamics using momentum to achieve better transient performance are also of interest. To overcome the potential scalability limitations of deterministic algorithms, future research will also explore hybrid optimization dynamics that incorporate stochastic phenomena via jumps.
%
\bibliographystyle{elsarticle-num}        
\bibliography{references.bib} 

\end{document}